\documentclass[12pt,journal]{new-aiaa}
\usepackage[utf8]{inputenc}
\usepackage{textcomp}
\usepackage{ulem} 
\usepackage{graphicx}
\usepackage{amsmath}
\usepackage[version=4]{mhchem}
\usepackage{siunitx}
\usepackage{longtable,tabularx}
\usepackage{color}
\usepackage{mathtools}
\usepackage{float}
\usepackage{caption,subcaption}
\usepackage{tikz}
\usepackage[many]{tcolorbox}
\usepackage{color}
\usepackage{epstopdf}
\usepackage[numbers,sort&compress]{natbib}
\usepackage{algorithm}
\usepackage{algorithmic}
\usepackage{setspace}
\usepackage[margin=1in]{geometry}
\def\1{\boldsymbol{1}}

\newtheorem{proposition}{Proposition}

\title{A Physics-Informed Indirect Method for Trajectory Optimization}

\author[1]{Kun Wang\footnote{PhD student, School of Aeronautics and Astronautics, Student member.}}
\author[1]{Fangmin Lu\footnote{PhD student, School of Aeronautics and Astronautics, Student member.}}
\author[1,2*]{Zheng Chen\footnote{Researcher, School of Aeronautics and Astronautics. email: \underline{z-chen@zju.edu.cn} (Corresponding author). }}
\author[1,2]{Jun Li\footnote{Professor, School of Aeronautics and Astronautics.}}
\affil[1]{School of Aeronautics and Astronautics, Zhejiang University, Hangzhou 310027, Zhejiang, China}
\affil[2]{Huanjiang Lab, Zhuji, 311816, Zhejiang, China}
\begin{document}

\maketitle

\begin{abstract}
The indirect method has long been favored for solving trajectory optimization problems
due to its ability to reveal the structure of the corresponding optimal control. However, providing an appropriate initial guess of the co-state vector remains challenging due to its lack of physical significance. To address this issue, this paper introduces a Physics-Informed Indirect Method (PIIM), \textcolor{blue}{for which all shooting variables are constrained into a narrow space by utilizing their physics information.}
Starting from the Time-Optimal Soft Landing Problem (TOSLP), an analytical estimation method for the minimum flight time is provided; we show that the co-state vector at the final time can be constrained on a unit 3-D hypersphere by eliminating the mass co-state and the numerical factor used in the initial co-state vector normalization; and the physical significance of the optimal control at the final time is exploited, allowing to further narrow down the solution space for the co-state vector to a unit 3-D octant sphere. Then, by incorporating the above physics-informed information into a shooting function that involves of propagating the dynamics of both states and co-states backward, the PIIM achieves faster and more robust convergence for the TOSLP compared to existing indirect methods.  In addition, the combination of the PIIM with a homotopy approach is proposed, allowing one to solve the \textcolor{blue}{fuel-optimal soft landing problem robustly. Finally, numerical simulations are presented to demonstrate and validate the developments of the paper.}
% demonstrate that,  while the success rate of the conventional indirect method for solving the FOSLP is $89.35\%$,  the proposed method takes a shorter time to find the solution to the FOSLP with a success rate of up to $100\%$.}
\end{abstract}
% \begin{IEEEkeywords}
% Indirect shooting method, Trajectory optimization, Initial co-state vector normalization, \textcolor{black}{Physical significance}, Backward propagation
% \end{IEEEkeywords}
\section{INTRODUCTION}
\label{intro}
A trajectory optimization problem aims to determine a sequence of inputs to a  dynamical system, while satisfying specific constraints and minimizing a predefined cost function. Due to the inherent nonlinearity of many real-world applications, numerical methods are commonly employed to generate optimal trajectories. The growing interest in space exploration, coupled with advancements in digital computing, has spurred academic and industrial sectors to develop efficient numerical methods to tackle intricate technical challenges. 

The first thorough classification of numerical methods for trajectory optimization was 
comprehensively conducted by Betts \cite{betts1998survey}, wherein direct and indirect methods were considered. The direct method transforms the trajectory optimization problem into a nonlinear programming  problem via collocation methods, which is then solved by the interior-point or sequential quadratic programming methods. However, this approach often entails a large number of variables, resulting in computationally expensive solutions.
On the other hand, the indirect method transforms the trajectory optimization problem into a Two-Point Boundary-Value Problem (TPBVP) based on the necessary conditions from Pontryagin’s Minimum Principle (PMP). The resulting TPBVP is usually solved by the indirect shooting method based on Newton’s method. Compared to the direct method, the indirect method offers two advantages. First, it reveals the structure of the optimal control through the necessary conditions. Second, it involves a smaller number of variables during the solution process, leading to faster convergence when the initial guess is sufficiently close to the optimal solution. However, the co-state vector is typically abstract with no physical significance, which makes providing an appropriate initial guess intricate and non-intuitive. Meanwhile, the indirect shooting method, when attempting to satisfy the terminal conditions, is very sensitive to the initial guess of the co-state vector \cite{conway2012survey}. Consequently, 
\textcolor{blue}{the convergence region for the initial co-state vector is often extremely narrow.}
% the region around the true co-state vector that leads to convergence to the optimal solution is often extremely narrow. 
Albeit not suitable for onboard implementation because the resulting solution is open-loop and solved offline, the indirect method has been widely used in obtaining the nominal trajectory, crucial for generating the real-time optimal solution via neural networks 
%\cite{cheng2018real} 
\cite{wang2024real}
or polynomial map guidance \cite{evans2024high}.
Therefore, extensive efforts have been made to find appropriate initial guesses that facilitate convergence and enhance the robustness of the indirect method.

Dixon and Biggs \cite{dixon1972advantages} introduced the adjoint control transformation to estimate physical controls and their derivatives, thus replacing the need for the initial co-state vector. This technique significantly reduced the sensitivity for orbital transfer problems. Subsequently, this technique was extended to other applications \cite{kluever1997optimal,yan1999initial,ranieri2005optimization,russell2007primer}. 
Another popular approach involves directly obtaining the analytical or approximate solution for the co-state vector by solving a simplified problem. For example, in \cite{thorne1996approximate}, an analytic co-state vector was obtained by simplifying the dynamics and disregarding the effects of the central body and mass variation. 
In \cite{cerf2012continuation}, a flat Earth model with constant gravity was employed to initialize the fuel-optimal orbit transfer problem. Additionally, in \cite{bonalli2017analytical}, for an optimal control problem concerning an endoatmospheric launch vehicle, an analytical initialization was first obtained, followed by a continuation method on the system dynamics to find the optimal solution of the original system.
In \cite{yang2019fast}, an analytic solution to the co-state vector from a gravity-free energy-optimal control problem was utilized to find the fuel-optimal trajectory for asteroid landing. In \cite{yang2020fuel}, the Lambert solution with irregular gravity was used to approximate the co-state vector for a fuel-optimal descent trajectory planning problem. 
After obtaining the analytical solution for a simplified one-dimensional landing problem, an adaptive homotopy process was designed to link the one-dimensional problem with the three-dimensional one in \cite{cheng2022real}. By using some reduced trajectory optimization problems as transitions, some final states and the co-state vector were initialized in analytical forms for the original minimum-time low-thrust problem in \cite{wu2021minimum}. In a low-thrust minimum-time station change problem \cite{zhao2016initial}, the tangential-thrust control was considered as the optimal control candidate, and the final time and co-state vector were approximated. 
Under an assumption of constant gravitational acceleration, the concept of zero-effort-miss/zero-effort-velocity has been used to derive closed-form suboptimal solutions \cite{guo2013applications,wang2021two}. 
In \cite{zhang2014fuel},  based on analytical results of a minimum-time problem, the optimal co-states were estimated to facilitate the solving procedure for the fuel-optimal lunar ascent problem. Recently, in virtue of the Hamiltonian being a linear function of the co-state vector,  the optimality conditions were transformed into linear forms of the co-state vector, and the  initial co-state vector was estimated by solving a set of linear algebraic equations \cite{wu2022analyticalcostate}.
In addition, shape-based methods, describing the geometric shape of the trajectory using mathematical expressions with tunable parameters, have also been widely used for trajectory optimization problems \cite{lunghi2015semi,taheri2017co,jiang2017improving,huo2022fast}. Specifically, shape-based methods have been applied to estimate the co-state vector for low-thrust transfer problems \cite{taheri2017co,jiang2017improving,huo2022fast}. Moreover, heuristic techniques have been used to initiate the co-state vector \cite{pontani2014optimal,hecht2023particle}.

In contrast to the approximation-based methods in the preceding paragraph, some works focus on reducing the solution space for the co-state vector. Lu et al. \rm{\cite{lu2008rapid}} analyzed a transversality condition in the optimal ascent problem and revealed that the co-state vector could be scaled by any positive constant without altering the optimal trajectory. Then, the Initial Co-state Vector Normalization (ICVN) was proposed in \cite{jiang2012practical}. In this work, through multiplying the cost function by a positive unknown numerical factor, the trajectory optimization problem was homogenized w.r.t. the co-state vector, effectively constraining the initially unbounded co-state vector to lie on a unit hypersphere. While originally designed for the fuel-optimal low-thrust problem, the ICVN has been applied to time-optimal problems \cite{wijayatunga2023exploiting,guo2023minimum}.

The conventional indirect method propagates both states and co-states forward in time, which usually makes the shooting function very sensitive to the initial guess of the co-state vector. Therefore, some works have aimed to mitigate this sensitivity. For instance, \cite{perez2018fuel} suggested that backward propagation could ensure the fulfillment of terminal conditions \textcolor{blue}{more easily compared to the forward propagation}. Meanwhile, some assumptions were made at the final time to reduce the solution space. \textcolor{blue}{In \cite{sidhoum2023indirect}, a forward-backward method was proposed to divide the trajectory by an intermediate point into two segments, and the trajectory optimization problem was formulated as a multi-PBVP rather than a TPBVP. During the solving process, the trajectory before the intermediate point was propagated forward, but the trajectory after the intermediate point was propagated backward. This forward-backward process allows for reducing the sensitivities of the shooting function.}

% As summarized, without simplifying the dynamics, the ICVN is a powerful technique to significantly reduce the solution space for the initial co-state vector down to a unit hypersphere. However, constraining the initial co-state vector to stay on the unit hypersphere may not be sufficient to ensure a very high convergence rate, as will be shown in Section \ref{NNNNNN}. In addition, backward propagation of the dynamics can be used to reduce the sensitivity of the shooting function. Nevertheless, the advantages of backward propagating the dynamics have not been well studied. 
\textcolor{blue}{
In summary, without simplifying the dynamics, the ICVN stands out as a powerful method for significantly reducing the solution space for the initial co-state vector down to a unit hypersphere. However, confining the initial co-state vector to remain on the unit hypersphere may not always guarantee a very high convergence rate, as will be shown in Section \ref{NNNNNN}. Additionally, backward propagation of the dynamics is able to facilitate convergence of the shooting function because the terminal conditions can be fulfilled more easily. Despite this potential, the benefits of backward propagating the dynamics have not been extensively explored.}

On the other hand, continuation or homotopy methods have been frequently embedded into the indirect method to facilitate convergence \cite{cerf2012continuation, bonalli2017analytical, yang2019fast, jiang2012practical, wijayatunga2023exploiting, sidhoum2023indirect, tang2018fuel, li2021homotopy, wang2022indirect}. When solving the fuel-optimal problem, a more straightforward subproblem, usually the energy-optimal problem, is often formulated within the continuation framework \cite{jiang2012practical, tang2018fuel, guo2024multi}. However, the significance of the energy-optimal problem is probably limited to initializing the fuel-optimal problem because the thrust magnitude in the energy-optimal problem is continuously changing, potentially rendering it infeasible for spacecraft operations. In contrast, the time-optimal solution provides the minimum flight time for any feasible trajectory \cite{yang2017rapid}.
Therefore, it is natural to investigate the direct connection between the time- and fuel-optimal problems.

Due to these reasons, in this paper we propose a Physics-Informed Indirect Method (PIIM) to efficiently and robustly solve two commonly studied trajectory optimization problems: the Time-Optimal Soft Landing Problem (TOSLP) and the Fuel-Optimal Soft Landing Problem (FOSLP). 
\textcolor{blue}{The key of the PIIM is that it provides all shooting variables, including the minimum flight time and unknown co-states, with physics-informed information. For the TOSLP,  unlike the approach to generate the minimum flight time iteratively by combining extrapolation and bisection methods \cite{yang2017rapid}, we provide an accurate estimation of the minimum flight time. We then show that the mass co-state and the numerical factor used in the ICVN can be eliminated by analyzing the dynamics at the final time. In contrast to existing works \cite{perez2018fuel,sidhoum2023indirect}, 
}one of the highlights of the PIIM is that, by propagating the dynamics of both states and co-states backward, it leverages the physical significance of the co-state vector at the final time to narrow down its solution space.
In this way, \textcolor{blue}{all the physics-informed information is embedded into a shooting function formulated within the framework of the PIIM. Consequently, the solution space for all shooting variables} can be reduced, allowing to facilitate convergence and improve the robustness of the indirect method. We summarize the main contributions of the paper as follows:
\begin{itemize}
  \item[1)] \textcolor{black}{ 
    We show that the co-state vector at the final time can be constrained on a unit 3-D hypersphere by eliminating the mass co-state and the numerical factor for the TOSLP.}
  \item[2)] \textcolor{black}{
 The physical significance of the optimal control at the final time is exploited, allowing to narrow down the solution space for the co-state vector at the final time to a unit 3-D octant sphere. Then, the PIIM can be implemented by incorporating the above physics-informed information into a shooting function that involves of propagating the dynamics of both states and co-states backward.}
  \item[3)] \textcolor{black}{We find that the final time tends to converge to a negative value if the initial guess of the co-state vector is not accurate enough.  Thus, a simple remedy strategy is proposed to guarantee that the final time remains positive in the shooting procedure.} 
  \item[4)] \textcolor{black}{A direct connection between the TOSLP and the FOSLP is established by combining the PIIM-based TOSLP with a homotopy approach, allowing to solve the FOSLP more quickly and robustly. }
  \end{itemize}

The structure of this paper is as follows. Section \ref{Pre} presents statements for the TOSLP and FOSLP. 
\textcolor{blue}{In Section \ref{SE:indirect_TOSLP}, the conventional indirect method based on the ICVN is briefly introduced. The PIIM for the TOSLP is detailed in Section \ref{PIIMTOSLP}.}
Section \ref{HomotopicApproach} outlines the conventional approach and the homotopy approach proposed to solve the FOSLP. In Section \ref{NNNNNN}, we present numerical simulations to demonstrate benefits of the proposed method.  Section \ref{Colus} concludes the paper.
\section{PROBLEM FORMULATION}
\label{Pre}
\begin{figure}[htb]
\begin{center}
\includegraphics[width=0.4\columnwidth]{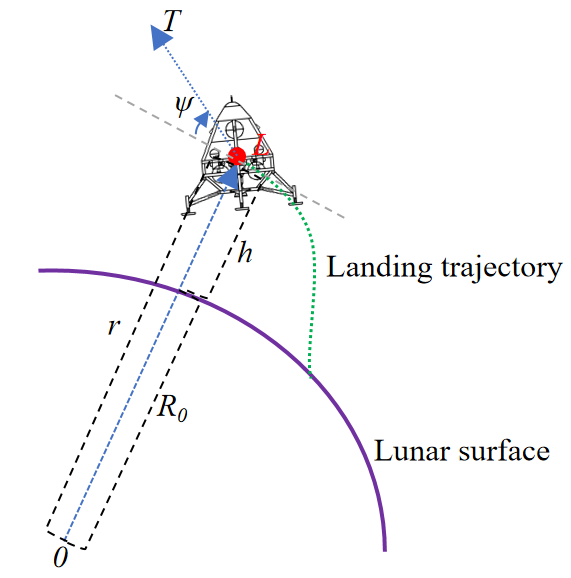}
\caption{\textcolor{black}{Coordinate system for the lunar soft landing.}}\label{Fig:frame}
\end{center}
\end{figure}

Consider the planar motion of a lunar lander described in Fig.~\ref{Fig:frame}. 
Assume that the Moon is a regular spherical body and the influences of the Moon’s rotation can be neglected. 
The origin $O$ is fixed at the Moon's center.
Denote by $r \ge R_0$ ($R_0$ is the radius of the Moon) the radial distance between the origin $O$ and the lunar lander $L$, thus the altitude of the lunar lander can be expressed as $h = r - R_0$. The lunar lander is propelled by an engine, whose thrust magnitude and thrust steering angle are adjustable. The thrust magnitude is denoted by $T \in [0,\textcolor{black}{T_\mathrm{m}}]$, where the constant $\textcolor{black}{T_\mathrm{m}}$ means the maximum thrust magnitude. Let $u\in[0,1]$ denote the engine thrust ratio, then we have $T=u\textcolor{black}{T_\mathrm{m}}$. The thrust steering angle, denoted by $\psi\in [-\frac{\pi}{2}, \frac{\pi}{2}]$, is defined to be the angle from the local horizontal line to the thrust vector. It is worth mentioning that since the landing site is not fixed, the optimal trajectory is rotatable. \textcolor{black}{Therefore, the initial range angle of the lunar lander, i.e., its initial angle information in the coordinate system, does not influence the solution.}
Then, the point-mass dynamics of the lunar lander can be described by \cite{liu2008optimal}
\begin{align}
\begin{cases}
\dot{r}(t) =  v(t),\\
\dot{v}(t) =  %\frac{u(t)\textcolor{black}{T_\mathrm{m}}}{m(t)}\sin\psi(t)-\frac{\mu}{r^2(t)} + r(t)\omega^2(t),\\
u(t)\textcolor{black}{T_\mathrm{m}}\sin\psi(t)/m(t)-\mu/{r^2(t)} + r(t)\omega^2(t),\\
%\dot{\theta}(t) = \omega(t), \\
\dot{\omega}(t) = -[u(t)\textcolor{black}{T_\mathrm{m}}\cos\psi(t)/m(t) + 2v(t)\omega(t)]/r(t),\\
\dot{m}(t) = %-\frac{u(t)\textcolor{black}{T_\mathrm{m}}}{\textcolor{black}{\textcolor{black}{I_\mathrm{sp}}}g_\mathrm{e}},
-u(t)\textcolor{black}{T_\mathrm{m}}/({\textcolor{black}{\textcolor{black}{I_\mathrm{sp}}}g_\mathrm{e}}),
\end{cases}
\label{LunarLander:DiffEqution}
\end{align}
where $t \ge 0$ is the time, the \textcolor{black}{over} dot denotes the differentiation w.r.t. time, $v$ is the speed along the direction of $OL$, i.e., the radial speed, and $\omega$ represents the angular velocity. Thus, the transverse speed is given by $\omega r$. $\mu$ represents the gravitational constant of the Moon, and $m$ is the mass of the lunar lander. The constant $\textcolor{black}{I_\mathrm{sp}}$ denotes the specific impulse of the lunar lander's \textcolor{black}{propulsion system}, while $g_\mathrm{e}$ represents the Earth's gravitational acceleration at sea level. The initial condition of the lunar lander at time $t_0 = 0$ is specified by
\begin{align}
r(0) = r_0, ~v(0) = v_0, ~\omega(0) = \omega_0, ~m(0) = m_0.
\label{InitialState}
\end{align}

It is evident that at touchdown, we have $r(t_f) = R_0$ ($t_f$ is the free final time). Meanwhile, we expect to have zero speed at touchdown, namely $v(t_f) = 0$ and $\omega(t_f) = 0$. Hence, the terminal condition can be expressed as
\begin{align}
r(t_f) = R_0, ~v(t_f) = 0, ~\omega(t_f) = 0.
\label{FinalState}
\end{align}

To improve numerical conditioning, we use $R_0$, $\sqrt\frac{\mu}{R_0}$, $m_0$, \textcolor{black}{$\sqrt\frac{{R^3_0}}{\mu}$}, and \textcolor{black}{$\frac{m_0 \mu}{{R^2_0}}$} to normalize $r$, $v$, $m$, $t$, and $\textcolor{black}{T_\mathrm{m}}$, respectively. Note that  $\mu$ will be normalized to $1$. To avoid abuse of notation, we continue to use the same notation as in (\ref{LunarLander:DiffEqution}) for the dimensionless counterpart hereafter.

To address the TOSLP, it amounts to finding control variables $u(t)$, $\psi(t)$, and the final time $t_f$, to steer the dynamical system in (\ref{LunarLander:DiffEqution}) from the initial condition in (\ref{InitialState}) to the terminal condition in (\ref{FinalState}), such that the cost function, 
\begin{align}
J_T = \int_0^{t_{f}} 1~\mathrm{d}t,
\label{EQ:cost_time}
\end{align}
is minimized, where the subscript $T$ is used in the context of TOSLP when necessary.

Unlike the TOSLP, the FOSLP minimizes the following cost function:
\begin{align}
J_F = \int_0^{t_{f}} u~\mathrm{d}t,
\label{EQ:cost_fuel}
\end{align}
where the subscript $F$ is used in the context of FOSLP when necessary.
\section{\textcolor{blue}{CONVENTIONAL INDIRECT METHOD FOR THE TOSLP}}\label{SE:indirect_TOSLP}
Let  $\boldsymbol{x} = [r,v,\omega,m]\textcolor{black}{^T}$ be the state vector \textcolor{black}{and} $\boldsymbol{p}_x = [p_r, p_v, p_{\omega},p_m]\textcolor{black}{^T}$ the co-state vector. 
\textcolor{black}{The nonlinear dynamics} in (\ref{LunarLander:DiffEqution}) can be rewritten \textcolor{black}{as}
\begin{align}
\dot{\boldsymbol x}(t) = \boldsymbol f(\boldsymbol x,  u, \psi, t),
\label{EQ:rewrite}
\end{align}
where $\boldsymbol f:\mathbb{R}^4 \times \mathbb{R} \times \mathbb{R} \times \mathbb{R}_0^+ \rightarrow \mathbb{R}^4$ is the smooth vector field defined in (\ref{LunarLander:DiffEqution}). 

\textcolor{blue}{Denote by $p_0 > 0$ an arbitrary positive numerical factor. The cost function in (\ref{EQ:cost_time}) is changed to}
\begin{align}
J_{\text{T,ICVN}} = p_0~\int_0^{t_{f}} 1~\mathrm{d}t.
\label{EQ:new_time_index}
\end{align}
The Hamiltonian, denoted by $\mathscr{H}_{\text{T,ICVN}}$, is formulated as 
\begin{align}
    \mathscr{H}_{\text{T,ICVN}} = 
    p_r v + p_v (\frac{u\textcolor{black}{T_\mathrm{m}}}{m}\sin\psi-\frac{\mu}{r^2} + r\omega^2) 
    +  p_\omega[-(\frac{u\textcolor{black}{T_\mathrm{m}}}{m}\cos\psi + 2v\omega)/r] + p_m(-\frac{u\textcolor{black}{T_\mathrm{m}}}{\textcolor{black}{I_\mathrm{sp}}g_\mathrm{e}}) + p_0.        
\label{EQ:Ham_ICVN}
\end{align}
% The Hamiltonian is expressed \textcolor{black}{as}
% \begin{align}
% \begin{split}
% \mathscr{H}&=  
% p_r v + p_v (\frac{u\textcolor{black}{T_\mathrm{m}}}{m}\sin\psi-\frac{\mu}{r^2} + r\omega^2) \\
% &+ p_\omega[-(\frac{u\textcolor{black}{T_\mathrm{m}}}{m}\cos\psi + 2v\omega)/r] + p_m(-\frac{u\textcolor{black}{T_\mathrm{m}}}{\textcolor{black}{I_\mathrm{sp}}g_\mathrm{e}}) + 1. 
% \end{split}
% \label{EQ:Ham}
% \end{align}
According to PMP \cite{Pontryagin}, \textcolor{black}{we have}
\begin{align}
\begin{cases}
\dot{p}_r(t) =  -\frac{2p_v\mu}{r^3} -  p_v{\omega}^2 -p_{\omega}[(\frac{u\textcolor{black}{T_\mathrm{m}}}{m}\cos\psi + 2v\omega)/r^2],\\
\dot{p}_v(t) =   -p_r + 2p_{\omega}\omega/r,\\
%\dot{p}_{\theta}(t) = -\frac{\partial \mathscr{H}}{\partial \theta} = 0,\\
\dot{p}_{\omega}(t) =  -2 p_v r\omega 
+ 2 p_{\omega} v/{r}, \\
\dot{p}_m(t) = {p_v u \textcolor{black}{T_\mathrm{m}}\sin\psi}/{m^2}  - p_{\omega} u \textcolor{black}{T_\mathrm{m}\cos\psi}/{(m^2 r)} .
\end{cases}
\label{EQ:dot_p}
\end{align}
% The optimal thrust steering angle follows that
% \begin{align}
% \frac{\partial \mathscr{H}}{\partial \psi} = 0. \label{EQ:direction}
% \end{align}
% Explicitly rewriting (\ref{EQ:direction}) leads to
% \begin{align}
% \psi(t) =  \arctan [-\frac{p_v(t)r(t)}{p_\omega(t) }].
% \label{EQ:optimal_psi}
% \end{align}
Minimizing $\mathscr{H}_{\text{T,ICVN}}$ w.r.t. $\psi$  implies
\begin{align}
\left[ \begin{array}{c}
\sin\psi \\
\cos\psi
\end{array}
\right]  = - \frac{1}{\sqrt{{p^2_v}+(-\frac{p_w}{r})^2}}\left[ \begin{array}{c}
p_v \\
-\frac{p_w}{r}
\end{array}
\right].
\label{EQ:optimal_H1}
\end{align}
For the optimal engine thrust ratio, we have
\begin{align*}
u(t)= 
\left\{ 
    \begin{array}{lc}
        1,~S(t) \leq 0 \\
        0,~S(t) > 0\\
        \text{undetermined},~S(t) = 0\\
    \end{array}
\right.
%\label{EQ:magnitude}
\end{align*}
where $S(t)$ denotes the switching function, defined as
 \begin{align}
% S(t) = \frac{p_v(t) \textcolor{black}{T_\mathrm{m}}}{m(t)} \sin\psi(t) -\frac{p_\omega(t) \textcolor{black}{T_\mathrm{m}}}{m(t)r(t)}\cos\psi(t)-\frac{p_m(t)\textcolor{black}{T_\mathrm{m}}}{\textcolor{black}{I_\mathrm{sp}}g_\mathrm{e}}.
S(t) = \frac{\partial \mathscr{H}_{\text{T,ICVN}}}{\partial u}.
\label{EQ:SF}
\end{align}
Substituting (\ref{EQ:optimal_H1}) into (\ref{EQ:SF}) leads to
 \begin{align*}
S(t) = -\textcolor{black}{T_\mathrm{m}}\left\{\frac{1}{m(t)} \sqrt{{p^2_v(t)}+[{\frac{p_\omega(t)}{r(t)}}]^2} +\frac{p_m(t)}{\textcolor{black}{I_\mathrm{sp}}g_\mathrm{e}}\right\}.
%\label{EQ:SF_real}
\end{align*}
\textcolor{blue}{The free final mass of the lunar lander implies}
% As the final mass of the lunar lander is free, the transversality condition 
\begin{align}
p_{m}(t_f) = 0.
\label{EQ:Transversality1}
\end{align}
With some algebraic manipulations, it can be shown that, combining the dynamics for ${p}_m$, (\ref{EQ:optimal_H1}) with (\ref{EQ:Transversality1}) leads to $p_{m}(t) > 0, \forall~t \in [0,t_f).$ Therefore, we have that $S(t)$ is always negative and
thus:
\begin{align}
u(t)\equiv 1, \forall~t \in [0,t_f].
\label{EQ:magnitudenew}
\end{align}
Because the Hamiltonian in (\ref{EQ:Ham_ICVN}) does not contain time explicitly and the final time is free, it holds that
\begin{align}
  \mathscr{H}_{\text{T,ICVN}} (t)\equiv 0,\ \forall~ t\in[0,t_f].
\label{EQ:hamiszero}
\end{align}

Meanwhile, $p_0$ and the initial co-state vector are normalized to stay on a unit 5-D hypersphere, i.e., 
\begin{align}
p^2_0 + p^2_{r_0} + p^2_{v_0} + p^2_{\omega_0} + p^2_{m_0}= 1,
\label{EQ:Ham_ICVN_norm}
\end{align}
which is called the {\it normalization condition}.

Then, the initial co-state vector and final time can be found by solving the following shooting function:
\begin{align}
    \boldsymbol{\Phi}^f_T(\boldsymbol{z}^{f}_{\text{T,ICVN}}) = [r(t_f)-R_0; v(t_f); \omega(t_f); p_{m}(t_f);
    p^2_0 + p^2_{r_0} + p^2_{v_0} + p^2_{\omega_0} + p^2_{m_0}-1;\mathscr{H}_{\text{T,ICVN}}(t_f)]=\boldsymbol{0},    
\label{EQ:TPBVP_law_ICVN}
\end{align}
where $\boldsymbol{z}^{f}_{\text{T,ICVN}} = [p_{r_0}, p_{v_0}, p_{\omega_0}, p_{m_0}, p_0, t_f]\textcolor{black}{^T}$ is the shooting vector. The superscript $f$ implies that the dynamics are propagated {\it forward};
Notice that $p_0$ cannot be $1$, otherwise $[p_{r_0},p_{v_0},p_{\omega_0},p_{m_0}]\textcolor{black}{^T} = \boldsymbol{0}$, which contradicts PMP. 

% \begin{align}
% \boldsymbol{\Phi}^f_T(\boldsymbol{z}^f_T) = [r(t_f)-R_0; v(t_f); \omega(t_f); p_{m}(t_f); \mathscr{H}(t_f)]=\boldsymbol{0},
% \label{EQ:TPBVP_law}
% \end{align}
% where $\boldsymbol{\Phi}^f_T$ is the shooting function with the superscript $f$ implying that the dynamics of states and co-states are propagated {\it forward}; $\boldsymbol{z}^f_T = [p_{r_0}, p_{v_0},  p_{\omega_0},p_{m_0},t_f]\textcolor{black}{^T}$ is called the \textcolor{black}{shooting vector}, in which $[p_{r_0}, p_{v_0},  p_{\omega_0},p_{m_0}]\textcolor{black}{^T}$ is the initial guess of the co-state vector $\boldsymbol{p_x}$   
% and $t_f$ is the initial guess of the final time. 
\textcolor{blue}{Denote by $t_\mathrm{max}$ the maximum value for the final time $t_f$. Since the thrust magnitude remains at its maximum during landing, $t_\mathrm{max}$ can be given by \cite{lunghi2015semi} (for simplicity, we ignore the dry mass of the lunar lander)
\begin{align}
  t_\mathrm{max} = m_0 \frac{I_\mathrm{sp}g_\mathrm{0}}{T_\mathrm{m}}.
  \label{EQ:tf_max}
\end{align}}
% \textcolor{black}{Apparently}, the solution space for (\ref{EQ:TPBVP_law}) is given as
% \begin{align}
%   \begin{split}
%     &p_{r_0} \in (-\infty,+\infty), p_{v_0} \in (-\infty,+\infty), p_{\omega_0} \in (-\infty,+\infty), \\
%     &p_{m_0} \in (0,+\infty), t_f \in (0,\textcolor{blue}{t_\mathrm{max}]}.
%   \end{split}
% \label{EQ:solutionspace_simple}
% \end{align}

Thus, the solution space for (\ref{EQ:TPBVP_law_ICVN}) becomes
\begin{align}
    p_{r_0} \in (-1,1), p_{v_0} \in (-1,1), p_{\omega_0} \in (-1,1), 
    p_{m_0} \in (0,1), p_0 \in (0,1), t_f \in (0,\textcolor{blue}{t_\mathrm{max}]}.    
\label{EQ:solutionspace_ICVN}
\end{align}

% Solving (\ref{EQ:TPBVP_law}) could be extremely challenging due to the boundless nature of the first four variables in $\boldsymbol{z}^f_T$. 
% \textcolor{black}{
% The stationary condition becomes}
% \begin{align}
% \mathscr{H}_{\text{T,ICVN}}(t_f) = 0. 
% \label{EQ:Ham_ICVN_zero}
% \end{align}
% \rm{Therefore}, the shooting function in (\ref{EQ:TPBVP_law}) becomes
% \begin{align}
%   \begin{split}
%     &\boldsymbol{\Phi}^f_T(\boldsymbol{z}^{f}_{\text{T,ICVN}}) = [r(t_f)-R_0; v(t_f); \omega(t_f); p_{m}(t_f);\\
%     &p^2_0 + p^2_{r_0} + p^2_{v_0} + p^2_{\omega_0} + p^2_{m_0}-1;\mathscr{H}_{\text{T,ICVN}}(t_f)]=\boldsymbol{0},    
%   \end{split}
% \label{EQ:TPBVP_law_ICVN}
% \end{align}
% where $\boldsymbol{z}^{f}_{\text{T,ICVN}} = [p_{r_0}, p_{v_0}, p_{\omega_0}, p_{m_0}, p_0, t_f]\textcolor{black}{^T}$ is the \textcolor{black}{shooting vector}. 
% Notice that $p_0$ cannot be $1$, otherwise we have $[p_{r_0},p_{v_0},p_{\omega_0},p_{m_0}]\textcolor{black}{^T} = \boldsymbol{0}$ according to (\ref{EQ:Ham_ICVN_norm}). This contradicts PMP. Thus, the solution space for (\ref{EQ:TPBVP_law_ICVN}) becomes
% \begin{align}
%   \begin{split}
%     &p_{r_0} \in (-1,1), p_{v_0} \in (-1,1), p_{\omega_0} \in (-1,1), \\
%     &p_{m_0} \in (0,1), p_0 \in (0,1), t_f \in (0,\textcolor{blue}{t_\mathrm{max}]}.    
%   \end{split}
% \label{EQ:solutionspace_ICVN}
% \end{align}
%\vspace{-0.7cm}
\section{\textcolor{blue}{PIIM FOR THE TOSLP}}\label{PIIMTOSLP}
\textcolor{blue}{
In this section, we first provide an analytical estimation of the minimum flight time $t_f$. Then, we show that the mass co-state $p_m$ can be eliminated and the numerical factor $p_0$ does not change the optimal solution by using a physical fact at the final time. This not only ensures that the initial co-state vector is constrained on a unit 3-D hypersphere, but also leads to the consideration of propagating the dynamics backward. 
The physical significance of the optimal control at the final time is exploited, allowing to narrow down the solution space. As a result, all the physics-informed information brought by the PIIM is incorporated into the shooting function.
}
\subsection{\textcolor{blue}{Estimation of the minimum flight time}}\label{estimationOF_flighttime}
By taking the lunar surface as the zero of potential energy, the initial energy of the lunar lander, denoted by $E_0$, is
\begin{align*}
E_0 = \frac{1}{2}m_0[v^2_0 + (\omega_0 r_0)^2] + m_0\frac{\mu}{r^2_0}(r_0-R_0).
%\label{EQ:e0}
\end{align*}
It is clear that the energy of the lunar lander reaches zero at touchdown. Therefore, the energy variation during landing is $\Delta E = E_0 - 0= E_0.$ Denote by ${\Delta m}$ the fuel consumption during landing, and let ${\Delta \hat{m}}$ be the estimation of ${\Delta m}$. Since the thrust is mainly used to nullify the energy of the lunar lander, according to Tsiolkovsky's rocket equation \cite{burleson2002konstantin}, the fuel consumption can be estimated by 
\begin{align}
\Delta \hat{m} \approx  m_0 - m(t_f) \approx \eta m_0[1 - \exp ({-\frac{\Delta V}{\textcolor{black}{I_\mathrm{sp}}g_\mathrm{e}}})],
\label{EQ:delta_m}
\end{align}
where $\eta$ is a constant, and $\Delta V$ is the change in velocity, which can be approximated by
\begin{align*}
\frac{1}{2}m_0 \Delta^2 V \approx \Delta E.
%\label{EQ:delta_velocity}
\end{align*}
Due to the fact that the fuel consumption for finite thrust is slightly larger than that of impulsive  thrust regarding the same change in velocity $\Delta V$, we set $\eta=1.05$. 
% Therefore, the estimation of $m(\tau)|_{\tau=0}$, denoted by $\hat{m}(\tau)|_{\tau=0}$, can be obtained by
% \begin{align*}
%   \begin{split}
%     \hat{m}(\tau)|_{\tau=0} &= 1.05 m_0 \exp [{-\frac{\sqrt{v^2_0 + (\omega_0 r_0)^2 + 2\frac{\mu}{r^2_0}(r_0-R_0)}}{\textcolor{black}{I_\mathrm{sp}}g_\mathrm{e}}}]\\
%     &-0.05 m_0.    
%   \end{split}
% %\label{EQ:back_m0}
% \end{align*}
\textcolor{blue}{As a result, the estimation of the minimum flight time}
$t_f$, denoted by $\hat{t}_f$, can be given by
\begin{align}
t_f \approx
\hat{t}_f = \frac{\Delta \hat{m}}{|\dot{m}|} = \frac{\Delta \hat{m} \textcolor{black}{I_\mathrm{sp}}g_\mathrm{e}}{T_\mathrm{m}}.
\label{EQ:t_f}
\end{align}
\subsection{Elimination of the mass co-state and the numerical factor}\label{eliminated_pm0_p0}
Notice that $p_m$ is not included in the right-hand side of (\ref{LunarLander:DiffEqution}) and  (\ref{EQ:dot_p}). In addition, the optimal thrust steering angle $\psi$ in (\ref{EQ:optimal_H1}) is also independent of $p_m$. Meanwhile, the optimal engine thrust ratio remains at its maximum regardless of $p_m$. Hence, its initial value $p_{m_0}$ has no impact on the integration of the state equation in (\ref{LunarLander:DiffEqution}) and the co-state equation in (\ref{EQ:dot_p}) except for the dynamics of $p_{m}$ in (\ref{EQ:dot_p}). 

According to (\ref{EQ:dot_p}) and (\ref{EQ:optimal_H1}), the mass co-state variation, denoted by $\Delta p_m$, can be obtained by
\begin{equation}
\Delta p_m = p_m(t_f) - p_{m_0} = \int_0^{t_{f}} \dot{p}_m(t)~\mathrm{d}t 
= \int_0^{t_{f}} -\frac{u \textcolor{black}{T_\mathrm{m}}}{m^2(t)} \left\{ \sqrt{{p^2_v(t)}+[{\frac{p_\omega(t)}{r(t)}}]^2}\right\} ~\mathrm{d}t.
\label{EQ:deltapm}
\end{equation}
Since all the elements on the right-hand side of (\ref{EQ:deltapm}) are independent of $p_{m_0}$, it is evident that $\Delta p_m$ is also independent of $p_{m_0}$. Consequently, once the terminal condition in (\ref{FinalState}) and the stationary condition in (\ref{EQ:hamiszero}) are met, $p_m(t_f)=0$ will automatically hold true by selecting $p_{m_0} = -\Delta p_m$. 
It is important to recall \textcolor{black}{that} $p_0$ can be any positive number without changing the nature of the optimal solution. Once the optimal solution is obtained, we set a new numerical factor, denoted by $p_0'$, via multiplying $p_0$ by a scaling factor \textcolor{black}{$k = \frac{1}{\sqrt{1-(\Delta p_m)^{2}}} > 1$}, i.e.,
\begin{equation*}
p_0' =  k p_0 =  \frac{1}{\sqrt{1-(\Delta p_m)^{2}}} p_0.
%\label{EQ:p_0_ICVN_p0}
\end{equation*}
Since the dynamics of states and co-states, as well as the optimal control equations are homogeneous to $p_0$ and the initial co-state vector, we have
\begin{equation*}
p_{r_0}' =  k p_{r_0}, p_{v_0}' =  k p_{v_0}, p_{\omega_0}' =  k p_{\omega_0},
%\label{EQ:p_0_ICVN_newco-states}
\end{equation*}
where the prime indicates the solution related to the new numerical factor $p_0'$.
Meanwhile, along the same optimal trajectory with different numerical factors $p_0$ and $p_0'$, the following relation holds: 
\begin{align}
    p_{r}'(t_f') =  k p_{r}(t_f), p_{v}'(t_f') =  k p_{v}(t_f), 
    p_{\omega}'(t_f') =  k p_{\omega}(t_f), t_f' = t_f, m'(t_f') = m(t_f).    
\label{EQ:p_0_ICVN_newstates}
\end{align}
By substituting (\ref{FinalState}), (\ref{EQ:optimal_H1}), and (\ref{EQ:Transversality1}),  (\ref{EQ:magnitudenew}) into (\ref{EQ:Ham_ICVN}), and considering (\ref{EQ:hamiszero}), we can deduce
\begin{align}
    \mathscr{H}_{\text{T,ICVN}}(t_f) = -\frac{\textcolor{black}{T_\mathrm{m}}}{m(t_f)}\sqrt{p^2_v(t_f) + p^2_{\omega}(t_f)} 
    - p_v(t_f) +  p_0 = 0.    
\label{EQ:p_0_ICVN}
\end{align}
Combining (\ref{EQ:p_0_ICVN_newstates}) with (\ref{EQ:p_0_ICVN}), we immediately have 
\begin{equation*}
 \mathscr{H}'_{\text{T,ICVN}}(t_f')= 
k\left\{-\frac{\textcolor{black}{T_\mathrm{m}}}{m(t_f)}\sqrt{p^2_v(t_f) + p^2_{\omega}(t_f)} - p_v(t_f) +  p_0\right\},
%\mathscr{H}'_{\text{T,ICVN}}(t_f') = k\mathscr{H}_{\text{T,ICVN}}(t_f),
%\label{EQ:p_0_ICVN_new}
\end{equation*}
which will be equal to zero once (\ref{EQ:p_0_ICVN}) is met. Additionally, in view of (\ref{EQ:Ham_ICVN_norm}), the scaled numerical factor and initial co-state vector satisfy 
\begin{align*}
p'^2_0 + p'^2_{r_0} + p'^2_{v_0} + p'^2_{\omega_0} = k^2 [1-(\Delta p_m)^{2}] = 1.
%\label{EQ:p_0_pm0}
\end{align*}

Hence, for any $p_0$, we can always set a new numerical factor $p_0' = k p_0$ such that $p_{m_0}$ can be eliminated from the {\it normalization condition} in (\ref{EQ:Ham_ICVN_norm}). To avoid notation abuse, we will still use the same notation for the scaled counterpart. Consequently,  (\ref{EQ:Ham_ICVN_norm}) becomes
\begin{align}
p^2_0 + p^2_{r_0} + p^2_{v_0} + p^2_{\omega_0} = 1. 
\label{EQ:Ham_ICVN_norm_pm0}
\end{align}
Thus, (\ref{EQ:TPBVP_law_ICVN}) becomes
\begin{align}
    \boldsymbol{\Phi}^f_T(\boldsymbol{z}'^f_{\text{T,ICVN}}) = [r(t_f)-R_0; v(t_f); \omega(t_f); 
    p^2_0 + p^2_{r_0} + p^2_{v_0} + p^2_{\omega_0} -1; \mathscr{H}_{\text{T,ICVN}}(t_f)]=\boldsymbol{0},  
\label{EQ:TPBVP_law_ICVN_pm}
\end{align}
where $\boldsymbol{z}'^f_{\text{T,ICVN}} = [p_{r_0}, p_{v_0}, p_{\omega_0}, p_0, t_f]\textcolor{black}{^T}$ is the \textcolor{black}{shooting vector} independent of $p_{m_0}$.

Next, we shall show that $p_0$ does not influence the optimal solution. By \textcolor{blue}{using a physical fact,} \textcolor{black}{we show that} $p_0 > 0$ always holds true as long as (\ref{EQ:p_0_ICVN}) is satisfied.
\begin{proposition}\label{theorem_1}
For any $p_v(t_f)$ and $p_{\omega}(t_f)$, (\ref{EQ:p_0_ICVN}) always leads to a positive numerical factor $p_0$.
\end{proposition}
Proof. First, in order for (\ref{EQ:p_0_ICVN}) to hold true, it is necessary that $p_v(t_f)$ and $p_{\omega}(t_f)$ cannot both be zero. During the final phase of the landing, the thrust magnitude must be greater than the gravitational force to nullify the touchdown speed. Meanwhile, the gravitational acceleration at touchdown, i.e., $\frac{\mu}{r^2(t_f)}$ with $r(t_f)=R_0$,  is normalized to 1. Therefore, we have $T_\mathrm{m} > m(t_f)$. As a result, it is evident that
\begin{align*}
p_0 = \frac{\textcolor{black}{T_\mathrm{m}}}{m(t_f)}\sqrt{p^2_v(t_f) + p^2_{\omega}(t_f)} + p_v(t_f)  > 0
\end{align*}
always holds true, thereby completing the proof. $\square$

As stated in \cite{guo2023minimum},  $p_0$ can be set to a fixed value in advance and remains unchanged during the solution process. Therefore, according to Proposition \ref{theorem_1} and \cite{guo2023minimum},  (\ref{EQ:Ham_ICVN_norm_pm0}) does not need to be satisfied, and  (\ref{EQ:hamiszero}) can be replaced by the following equation:
\begin{align}
p^2_{r_0} + p^2_{v_0} + p^2_{\omega_0} = 1, 
\label{EQ:Ham_ICVN_norm_pm0_p0}
\end{align}
which is called the Simplified ICVN (SICVN), as done in \cite{guo2023minimum}. Then,  (\ref{EQ:TPBVP_law_ICVN_pm}) becomes 
\begin{align}
    \boldsymbol{\Phi}^f_T(\boldsymbol{z}^f_{\text{T,SICVN}}) = [r(t_f)-R_0; v(t_f); \omega(t_f);
    p^2_{r_0} + p^2_{v_0} + p^2_{\omega_0} -1]=\boldsymbol{0},    
\label{EQ:TPBVP_law_ICVN_pm_p0}
\end{align}
in which $\boldsymbol{z}^f_{\text{T,SICVN}} = [p_{r_0}, p_{v_0}, p_{\omega_0}, t_f]\textcolor{black}{^T}$ is the \textcolor{black}{shooting vector} independent of $p_{m_0}$ and $p_0$. \textcolor{black}{As a result}, $p_0$ can be determined by (\ref{EQ:p_0_ICVN}) once  (\ref{EQ:TPBVP_law_ICVN_pm_p0}) is solved. Moreover, the solution space for the co-state vector in (\ref{EQ:TPBVP_law_ICVN_pm_p0}) \textcolor{black}{is} 
\begin{align}
p_{r_0} \in (-1,1), p_{v_0} \in (-1,1), p_{\omega_0} \in (-1,1).
\label{EQ:solutionspace_ICVN_simple}
\end{align}
%\vspace{-0.6cm}
%\section{PIIM FOR THE TOSLP}\label{PIIM_section_toslp}
\subsection{Backward propagating of dynamics}\label{back_ward_propagate}
\textcolor{blue}{Since the developments of the preceding subsection are derived by analyzing the dynamics at the final time, we consider propagating the dynamics backward.}
% Notice that the dynamics of states and co-states in (\ref{EQ:rewrite}) and (\ref{EQ:dot_p}) are both propagating forward in time from $t = 0$ to $t= t_f$. 

Define a new time variable $\tau$ as below
\begin{align*}
\tau = t_f-t, t \in [0,t_f].
%\label{EQ:time_back}
\end{align*}
It is palpable that the dynamics can also be propagated backward in time, i.e., from $t=t_f$ to $t=0$. In such case, the dynamics are transformed into 
\begin{align}
\begin{cases}
\dot{\boldsymbol x}(\tau) = -\boldsymbol f(\boldsymbol x, u, \psi, \tau), \\
\dot{\boldsymbol{p}}_{x}(\tau)= \frac{\partial \mathscr{H}(\tau)}{\partial \boldsymbol x(\tau)}.
\end{cases}
\label{EQ:dynamics_backward}
\end{align}
For the TOSLP with the same boundary conditions, the initial condition in (\ref{InitialState}) becomes  the terminal condition when the dynamics are propagated backward, i.e.,
\begin{align}
    r(\tau) |_{\tau=t_f}= r_0,~v(\tau) |_{\tau=t_f} =v_0,
    %~\theta(\tau) |_{\tau=t_f} = \theta_0,
    \omega(\tau) |_{\tau=t_f} = \omega_0,~m(\tau) |_{\tau=t_f} = m_0.        
\label{EQ:dynamics_initial_terminal_back}
\end{align}
Likewise, the terminal condition in (\ref{FinalState}) becomes the \textcolor{black}{initial condition}, i.e.,
\begin{align}
r(\tau) |_{\tau=0} = R_0,~v(\tau) |_{\tau=0} = 0, ~\omega(\tau) |_{\tau=0} = 0. 
\label{EQ:dynamics_initial_back}
\end{align}
In such case, the variable $\tau$ represents the {\it time of flight}. Note that the optimal control equations for the dynamical system in (\ref{EQ:dynamics_backward}) remain the same as in (\ref{EQ:optimal_H1}) and (\ref{EQ:magnitudenew}) once the variable $t$ is replaced by $\tau$.

To propagate the dynamics in (\ref{EQ:dynamics_backward}) from $\tau=0$ to $\tau=t_f$, we need the unknown initial state $m(\tau) |_{\tau=0}$
and initial co-state vector $\boldsymbol{p}_{x}(\tau) |_{\tau=0} =[p_{r}(\tau), p_{v}(\tau), p_{\omega}(\tau)]\textcolor{black}{^T}|_{\tau=0}$. Notably, 
$p_{m}(\tau)|_{\tau=0}$ is eliminated because it does not alter the propagation on the remaining states and co-states.
Given (\ref{EQ:dynamics_initial_back}), $m(\tau) |_{\tau=0}$ and $\boldsymbol{p}_{x}(\tau) |_{\tau=0}$ to be determined, the final states and co-states can be calculated through integrating (\ref{EQ:dynamics_backward}) along with the optimal control equations.
Notice that along the optimal trajectory, $\mathscr{H}_{\text{T,ICVN}}(t) \equiv 0$ holds true for $t \in [0, t_f]$. For (\ref{EQ:dynamics_backward}), we consider the Hamiltonian $\mathscr{H}_{\text{T,ICVN}}(\tau)$ at $\tau = 0$ for the stationary condition. Recall that (\ref{EQ:Ham_ICVN_norm_pm0_p0}) can replace the stationary condition. 
 Thus, the TOSLP is transformed into a new TPBVP, which seeks to find $m(\tau) |_{\tau=0}$,  $\boldsymbol{p}_{x}(\tau) |_{\tau=0}$, and  $t_f$ such that  (\ref{EQ:dynamics_initial_terminal_back}) and (\ref{EQ:dynamics_initial_back}) are satisfied, i.e.,
\begin{equation}
\boldsymbol{\Phi}^b_T(\boldsymbol z^b_{\text{T,SICVN}}) = [r(\tau)|_{\tau=t_f}-r_0; v(\tau)|_{\tau=t_f}-v_0;
\omega(\tau)|_{\tau=t_f}-\omega_0;
m(\tau)|_{\tau=t_f} - m_0;
p^2_r(\tau)+p^2_v(\tau)+p^2_\omega(\tau)|_{\tau=0} - 1]=\boldsymbol{0},
\label{EQ:TPBVP_new}
\end{equation}
where the superscript $b$ denotes the {\it backward} propagation, and the \textcolor{black}{shooting vector} $\boldsymbol z^b_{\text{T,SICVN}}$ is
\begin{align}
\boldsymbol z^b_{\text{T,SICVN}} = [p_{r}(\tau),p_{v}(\tau), p_{\omega}(\tau), m(\tau), t_f]\textcolor{black}{^T}|_{\tau=0}. 
\label{EQ:shooting_new_back}
\end{align}

Compared with  (\ref{EQ:TPBVP_law_ICVN_pm_p0}), $\boldsymbol z^b_{\text{T,SICVN}}$ has one more unknown variable $m(\tau) |_{\tau=0}$, which is related to the fuel consumption during landing. 
\textcolor{blue}{Fortunately, $m(\tau)|_{\tau=0}$ can be obtained by $m(\tau)|_{\tau=0} = m_0 - {\Delta m}$ according to (\ref{EQ:delta_m}). In the next subsection, we focus on reducing the solution space for $[p_{r}(\tau),p_{v}(\tau), p_{\omega}(\tau)]^T|_{\tau=0}$.}
% In the following, we first provide the analytical estimations of $m(\tau)|_{\tau=0}$ and the final time $t_f$ in order to realize fast and robust convergence when solving (\ref{EQ:TPBVP_new}). Subsequently, we exploit the physical significance of the optimal control \textcolor{black}{at the final time} to further narrow down the solution space for (\ref{EQ:shooting_new_back}).
\subsection{Reducing the solution space by exploiting the optimal control}\label{reducing_solution}
Notice that the solution space for the co-state vector in (\ref{EQ:shooting_new_back}) is constrained on a unit 3-D hypersphere. In what follows, we shall show that the solution space can be further reduced to a unit 3-D octant sphere. 
\begin{proposition}\label{theorem_2}
Let a short time interval $t \in [t_f-t_1,t_f)$ ($t_1$ is a small positive number) denote the final phase of the landing until touchdown. \textcolor{black}{It} holds true that $p_v(\tau) < 0$ and $p_\omega(\tau) > 0$ for the interval $\tau = t_f - t  \in (0,t_1]$.
\end{proposition}
Proof. For the interval $t \in [t_f-t_1,t_f)$, in view of (\ref{LunarLander:DiffEqution}), the radial distance $r(t)$ will monotonically decrease to $R_0$, indicating that
$\dot{r}{(t)} <0$, 
which further implies that $v{(t)} < 0$.
To nullify the radial speed at touchdown, it is necessary to have $\dot v{(t)} > 0$.
Moreover, since the  angular velocity $\omega (t)$ is very close to zero in such \textcolor{black}{an} interval, the second equation in (\ref{LunarLander:DiffEqution}) becomes 
\begin{align*}
\dot v{(t)} \approx  \frac{T_\mathrm{m}}{m(t)} \sin \psi(t) - \frac{\mu}{r^2(t)} > 0.
\end{align*}
Since $r(t) \approx 1$ for  $t \in [t_f-t_1,t_f)$ and the estimation of $m(t)$, denoted by $\hat{m}(t)$, can be obtained by (\ref{EQ:delta_m}), we have 
\begin{align}
\sin \psi(t) > \frac{\hat{m}(t)}{T_\mathrm{m}}.
\label{EQ:end_dotv_approx}
\end{align}
Analogously, to nullify the angular velocity at touchdown, $\omega (t)$ should gradually decrease to zero, yielding 
$\dot \omega{(t)} < 0$.
In view of the fourth equation in (\ref{LunarLander:DiffEqution}), we have 
\begin{align}
\cos \psi{(t)} > 0.
\label{EQ:end_psi}
\end{align}
Combining (\ref{EQ:end_dotv_approx}) with (\ref{EQ:end_psi}), it is obvious that 
\begin{align}
\arcsin{\frac{\hat{m}(t)}{T_\mathrm{m}}} < \psi(t)< \frac{\pi}{2}.
\label{EQ:end_psi_cos}
\end{align}
Using the optimal thrust steering angle in (\ref{EQ:optimal_H1}) infers that
$
p_v(t) < 0~~\text{and}~~p_w(t) > 0, t \in [t_f-t_1,t_f),$
indicating that 
\begin{align}
p_v(\tau) < 0~~\text{and}~~p_w(\tau) > 0, \tau \in (0,t_1],
\label{EQ:end_pv_pw}
\end{align}
which completes the proof. $\square$

Next, let us analyze the solution space for $p_r(\tau)|_{\tau=0}$. In view of (\ref{EQ:optimal_H1}), by differentiating $\sin \psi$ w.r.t $\tau$, along with some
algebraic manipulations, we obtain
\begin{align}
\dot{\psi}(\tau) \cos \psi(\tau) |_{\tau=0} = \frac{p_r(\tau)p^2_\omega(\tau)}{[p^2_v(\tau)+p^2_\omega(\tau)]^{\frac{3}{2}}}|_{\tau=0}.
\label{EQ:end_pr}
\end{align}
Combining the second equation in (\ref{EQ:optimal_H1}) with (\ref{EQ:end_pr}) leads to
\begin{align}
\dot{\psi}(\tau) |_{\tau=0} = \frac{p_r(\tau)p_\omega(\tau)}{p^2_v(\tau)+p^2_\omega(\tau)}|_{\tau=0}.
\label{EQ:end_pr_new}
\end{align}
Since the thrust steering angle defines the body attitude of the lunar lander \cite{lu2023propellant}, it is desirable to ensure that the thrust steering angle at touchdown is such that the lunar lander has a vertical attitude. If not, a terminal descent phase is usually needed to ensure a vertical landing. In such case, to ensure a smooth and gentle landing, it is reasonable to assume that the thrust steering angle gradually increases towards $90$ deg, or at least tends to do so until touchdown. In other words, we have $\dot{\psi}(t) > 0$ for $t \in [t_f-t_1,t_f]$, which is equivalent to $\dot{\psi}(\tau) > 0$ for $\tau \in [0,t_1]$. Therefore,
in view of (\ref{EQ:end_pv_pw}) and (\ref{EQ:end_pr_new}), we obtain
\begin{align}
p_r(\tau) |_{\tau=0} > 0.
\label{EQ:end_pr_positive}
\end{align}

So far, according to (\ref{EQ:end_pv_pw}) and (\ref{EQ:end_pr_positive}), the solution space has been further reduced to
\begin{align}
p_{r}(\tau)|_{\tau=0} \in (0,1), p_{v}(\tau)|_{\tau=0} \in (-1,0), p_{\omega}(\tau)|_{\tau=0} \in (0,1).
\label{EQ:shooting_new_back_initial_less}
\end{align}
In the next section, to demonstrate the developments of the PIIM, we shall first present the conventional indirect method to solve the FOSLP by adopting the ICVN. Then, we will combine the PIIM with a homotopy approach to directly connect the TOSLP and the FOSLP.
\section{PROCEDURES FOR SOLVING THE FOSLP} \label{HomotopicApproach}
\subsection{\textcolor{blue}{Conventional indirect method using the ICVN}}
By introducing an arbitrary positive numerical factor $p_{0_F}$, the cost function in (\ref{EQ:cost_fuel}) is changed to
\begin{align}
J_{\text{F,ICVN}} = p_{0_F}~\int_0^{t_{f}} u~\mathrm{d}t.
\label{EQ:new_fuel_index}
\end{align}
\textcolor{black}{
The Hamiltonian, denoted by} $\mathscr{H}_{\text{F,ICVN}}$, is reconstructed as 
\begin{align*}
    \mathscr{H}_{\text{F,ICVN}} = 
    p_r v + p_v (\frac{u\textcolor{black}{T_\mathrm{m}}}{m}\sin\psi-\frac{\mu}{r^2} + r\omega^2) 
    +  p_\omega[-(\frac{u\textcolor{black}{T_\mathrm{m}}}{m}\cos\psi + 2v\omega)/r] 
     + p_m(-\frac{u\textcolor{black}{T_\mathrm{m}}}{\textcolor{black}{I_\mathrm{sp}}g_\mathrm{e}}) + p_{0_F} u.         
%\label{EQ:Ham_ICVN_fuel}
\end{align*}
Equations (\ref{LunarLander:DiffEqution}), (\ref{EQ:dot_p}), and (\ref{EQ:optimal_H1}) still hold true, expect that the optimal engine thrust ratio $u_F(t)$ becomes
\begin{align*}
u_F(t)= 
\left\{ 
    \begin{array}{lc}
        1,~S_F(t) \leq 0 \\
        0,~S_F(t) > 0\\
    \end{array}
\right.
%\label{EQ:magnitude_h_new}
\end{align*}
where $S_F(t)$ is the switching function satisfying
\begin{align*}
S_F(t)= p_{0_F} -\textcolor{black}{T_\mathrm{m}}\left\{\frac{1}{m(t)} \sqrt{{p^2_v(t)}+[{\frac{p_\omega(t)}{r(t)}}]^2} +\frac{p_m(t)}{\textcolor{black}{I_\mathrm{sp}}g_\mathrm{e}}\right\}.
%\label{EQ:reho_h_fuel}
\end{align*}
Since $u_F(t)$ is bang-bang, which may result in numerical difficulties, \textcolor{black}{we adopt the smoothing technique from \cite{wang2023new} to approximate $u_F(t)$, i.e.,} 
\begin{align}
u_F(t) \approx  u_F(t,\delta)  = \frac{1}{2}(1-\frac{S_F(t)}{\sqrt{|S_F(t)|^2 + \delta}}),
\label{EQ:reho_h_smoothing_fuel}
\end{align}
in which $\delta$ is a small positive constant.

Meanwhile, \textcolor{black}{we have the following \it{normalization condition}}:
\begin{align}
p^2_{0_F} + p^2_{r_0} + p^2_{v_0} + p^2_{\omega_0} + p^2_{m_0}= 1. 
\label{EQ:Ham_ICVN_norm_fuel}
\end{align}
The shooting function is
\begin{align}
    \boldsymbol{\Phi}^f_F(\boldsymbol{z}^{f}_{\text{F,ICVN}}(\delta)) = [r(t_f)-R_0; v(t_f);\omega(t_f);p_{m}(t_f);  
     p^2_{0_F}+ p^2_{r_0}+ p^2_{v_0} + p^2_{\omega_0} + p^2_{m_0}-1;\mathscr{H}_{\text{F,ICVN}}(t_f)]=\boldsymbol{0},    
\label{EQ:TPBVP_law_ICVN_fuel}
\end{align}
where $\boldsymbol{z}^{f}_{\text{F,ICVN}}(\delta) = [p_{r_0}, p_{v_0}, p_{\omega_0}, p_{m_0}, p_{0_F}, t_f]\textcolor{black}{^T}$. \textcolor{blue}{Since the lower bound for the thrust magnitude is zero, according to (\ref{EQ:t_f})}, 
the solution space for (\ref{EQ:TPBVP_law_ICVN_fuel}) is
\begin{align}
    p_{r_0} \in (-1,1), p_{v_0} \in (-1,1), p_{\omega_0} \in (-1,1),
    p_{m_0} \in (0,1), p_{0_F} \in (0,1), t_f \in [\hat{t}_f,+\infty).    
\label{EQ:solutionspace_ICVN_fuel}
\end{align}
% \textcolor{blue}{
% \subsection{Eliminating the numerical factor}
% Note that the switching function in (\ref{EQ:reho_h_fuel}) is related to $p_m$. In such case, the mass co-state can not be eliminated. In virtue of (\ref{EQ:Ham_ICVN}) and (\ref{EQ:Ham_ICVN_fuel}), it is easy to see that both equations can be reduced to the form of (\ref{EQ:p_0_ICVN}) because the 
% optimal engine thrust ratio has to be 1, i.e, $u(t_f) = 1$, for both TOSLP and FOSLP. 
% Analogous to Proposition \ref{theorem_1}, it is evident that the numerical factor $p_{0_F}$ can be eliminated for the FOSLP. Similarly, we use SICVN to denote the resulting shooting function. Therefore, the shooting function becomes
% \begin{align}
%   \begin{split}
%     &\boldsymbol{\Phi}^f_F(\boldsymbol{z}^{f}_{\text{F,SICVN}}(\delta)) = [r(t_f)-R_0; v(t_f);\omega(t_f);p_{m}(t_f);\\
%     &p^2_{r_0} + p^2_{v_0} + p^2_{\omega_0} + p^2_{m_0}-1]=\boldsymbol{0},    
%   \end{split}
% \label{EQ:TPBVP_law_SICVN_fuel}
% \end{align}
% where $\boldsymbol{z}^{f}_{\text{F,SICVN}}(\delta) = [p_{r_0}, p_{v_0}, p_{\omega_0}, p_{m_0}, t_f]^T$ is the shooting vector. The solution space is
% \begin{align}
%   \begin{split}
%     &p_{r_0} \in (-1,1), p_{v_0} \in (-1,1), p_{\omega_0} \in (-1,1),\\
%     &p_{m_0} \in (0,1), t_f \in (0,+\infty).    
%   \end{split}
% \label{EQ:solutionspace_SICVN_fuel}
% \end{align}
\subsection{Homotopy approach}
\textcolor{black}{
Since the PIIM can reduce the solution space down to be sufficiently small for the TOSLP, it is natural to directly connect the TOSLP with the FOSLP within the homotopy framework.}
Denote by $\kappa$ a homotopy parameter that connects  (\ref{EQ:new_time_index}) with  (\ref{EQ:cost_fuel}), i.e.,
\begin{align}
J_h = \int_0^{t_{f}} p_0\kappa + (1-\kappa) u~\mathrm{d}t,
\label{EQ:performance_homotopy}
\end{align}
where the subscript $h$ denotes the context of {\it homotopy}. It can be observed that $J_h$ becomes $J_{\text{T,ICVN}} = \int_0^{t_{f}} p_0~\mathrm{d}t$ in (\ref{EQ:new_time_index}) if $\kappa = 1$. On the other hand, $J_h$ becomes $J_F = \int_0^{t_{f}} u~\mathrm{d}t$ in (\ref{EQ:cost_fuel}) if $\kappa = 0$. In this manner, the TOSLP serves as the seeding problem, whose solution can be found using the \textcolor{black}{PIIM} developed in Section \ref{PIIMTOSLP}. Then, by decreasing $\kappa$ from $1$ to $0$, the homotopy process can be completed by using the preceding convergent solution as the initial guess. Next, we shall present the corresponding shooting functions during the homotopy process.

Notice that the dynamical model and boundary conditions remain the same as in Section \ref{Pre}. Considering the cost function in (\ref{EQ:performance_homotopy}), the Hamiltonian \textcolor{black}{is}
\begin{align*}
    \mathscr{H}_h(\kappa) = p_0\kappa + (1-\kappa) u+
    p_r v + p_v (\frac{u\textcolor{black}{T_\mathrm{m}}}{m}\sin\psi-\frac{\mu}{r^2} 
    + r\omega^2)+ p_\omega[-(\frac{u\textcolor{black}{T_\mathrm{m}}}{m}\cos\psi + 2v\omega)/r]+ p_m(-\frac{u\textcolor{black}{T_\mathrm{m}}}{\textcolor{black}{I_\mathrm{sp}}g_\mathrm{e}}) .         
%\label{EQ:Ham_homo}
\end{align*}
Therefore, (\ref{EQ:dot_p}) and (\ref{EQ:optimal_H1}) still hold true. Since the optimal engine thrust ratio $u_h(\tau,\kappa)$ is linear w.r.t. the state dynamics and cost function, it is reasonable to assume that $u_h(\tau,\kappa)$ is bang-bang, i.e.,
\begin{align*}
u_h(\tau,\kappa)= 
\left\{ 
    \begin{array}{lc}
        1,~S_h(\tau,\kappa) \leq 0 \\
        0,~S_h(\tau,\kappa) > 0\\
    \end{array}
\right.
%\label{EQ:magnitude_h}
\end{align*}
where $S_h(\tau,\kappa)$ is the switching function satisfying
\begin{align*}
    S_h(\tau,\kappa) =  -\textcolor{black}{T_\mathrm{m}}\left\{\frac{1}{m(\tau)} \sqrt{{p^2_v(\tau)}+[{\frac{p_\omega(t)}{r(\tau)}}]^2} +\frac{p_m(\tau)}{\textcolor{black}{I_\mathrm{sp}}g_\mathrm{e}}\right\}
    +1 - \kappa.   
%\label{EQ:reho_h}
\end{align*}
Once again, the smoothing  technique from \cite{wang2023new} is used to approximate $u_h(\tau,\kappa)$, i.e.,
\begin{align}
u_h(\tau,\kappa) \approx  u_h(\tau,\kappa,\delta)  = \frac{1}{2}(1-\frac{S_h(\tau)}{\sqrt{|S_h(\tau)|^2 + \delta}}),
\label{EQ:reho_h_smoothing}
\end{align}
in which $\delta$ is the same constant as in (\ref{EQ:reho_h_smoothing_fuel}).

The corresponding TPBVP amounts to solving the following shooting function:
\begin{equation}
\boldsymbol{\Phi}^b_h(\boldsymbol z^b_{\text{h,SICVN}}(\kappa,\delta)) = [r(\tau,\kappa,\delta)-r_0;v(\tau,\kappa,\delta)-v_0;
\omega(\tau,\kappa,\delta) -\omega_0;m(\tau,\kappa,\delta) -m_0;\mathscr{H}_h(\tau,\kappa,\delta)]|_{\tau=t_f}=\boldsymbol{0},
\label{EQ:TPBVP_new_homo}
\end{equation}
where 
$\boldsymbol z^b_{\text{h,SICVN}}(\kappa,\delta) = [p_{r}(\tau,\kappa,\delta),p_{v}(\tau,\kappa,\delta), p_{\omega}(\tau,\kappa,\delta), \\
m(\tau,\kappa,\delta), t_f(\kappa,\delta)]^T|_{\tau=0}$. It is important to note that, unlike the TOSLP where $p_m$ is not needed, the dynamics of $p_m$ must be included here because the switching function is related to $p_m$. When propagating the dynamics of $p_m$, it is clear that its initial value is $p_m(\tau)|_{\tau=0} = 0$.

\begin{figure*}[htbp]      
  \centering
  \includegraphics[scale=0.38]{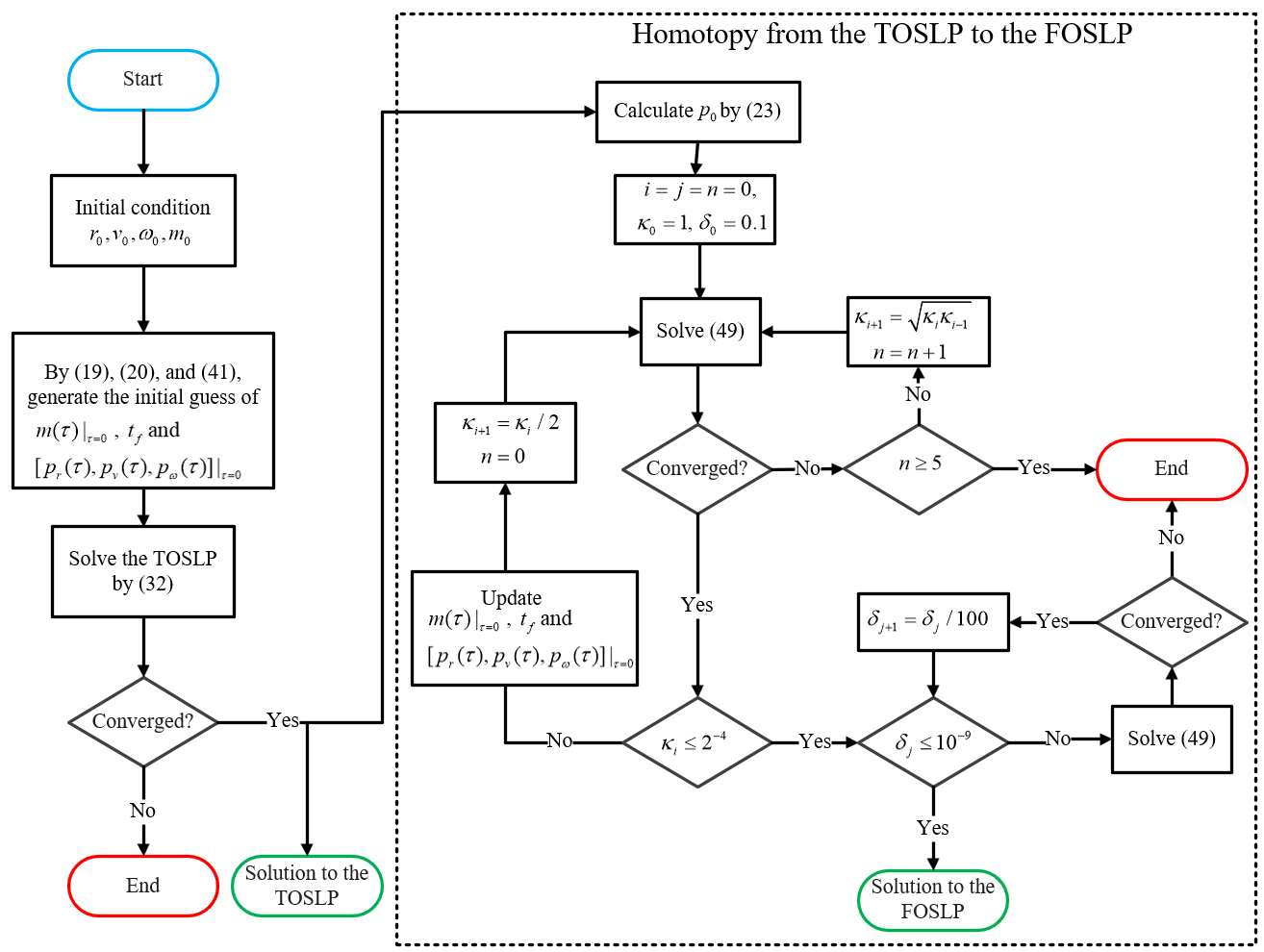}
  \caption{\textcolor{black}{Flowchart for solving the TOSLP and using its solution to solve the FOSLP via the homotopy process.}}\label{Fig:homo_flowchart}
  \end{figure*}

Keep in mind that (\ref{EQ:TPBVP_new}) is equivalent to (\ref{EQ:TPBVP_new_homo}) when $\kappa = 1$. In such case, the smoothing constant $\delta$ does not change the optimal thrust magnitude because it will remain at the maximum; meanwhile, the last equation in (\ref{EQ:TPBVP_new}) is used to replace the stationary condition, i.e., the last equation in (\ref{EQ:TPBVP_new_homo}). Instead of finding the solution to the FOSLP by directly solving (\ref{EQ:TPBVP_law_ICVN_fuel}), 
we shall present a homotopy process, as stated in Fig.~\ref{Fig:homo_flowchart},
to solve the FOSLP through starting with the \textcolor{black}{PIIM}-based TOSLP.  
In Fig.~\ref{Fig:homo_flowchart}, $i$ and $j$ are the updating indices for the homotopy parameters $\kappa$ and $\delta$, respectively. The parameter $n$ is the backtracking number for solving (\ref{EQ:TPBVP_new_homo}). This means that if the shooting function does not converge even after adopting the homotopy parameter updating procedure $5$ times, the algorithm fails.
\section{NUMERICAL SIMULATIONS}\label{NNNNNN}
% This section first presents some simulations on the TOSLP . Analytical estimations of fuel consumption and final time \textcolor{black}{required} for the \textcolor{black}{PIIM}-based shooting function in (\ref{EQ:TPBVP_new}) are verified. Then, three different methods are compared. \textcolor{black}{Lastly}, simulations on the FOSLP are used to showcase the developments of the proposed method. 
\textcolor{blue}{This section presents some simulations to showcase the developments in Sections  \ref{PIIMTOSLP} and \ref{HomotopicApproach}.
}
Before proceeding, we shall \textcolor{black}{define some constants}. The propulsion system of the lunar lander is specified by $\textcolor{black}{I_\mathrm{sp}} = 300$ s and $T_\mathrm{m} = 1,500$ N. The radius of the Moon is $R_0 = 1,738$ km, and $g_\mathrm{e}$ \textcolor{black}{is} $9.81$ $\rm m/s^2$. Additionally, the gravitational constant $\mu$ is $4.90275 \times 10 ^{12}$~$\rm m^3/s^2$. All the algorithms are implemented on a desktop equipped with an AMD EPYC 9684X 96-Core Processor @2.55 GHz and 128 GB of RAM. The absolute and relative tolerances for propagating the dynamics are set as $1.0 \times 10^{-9}$. \textcolor{blue}{The nonlinear equation solver {\it fsolve} is used to find the zero of the shooting function.}
The termination tolerance and maximum number of iterations for {\it fsolve} are set as $1.0 \times 10^{-9}$ and $300$, respectively.
\subsection{Simulations on the TOSLP}
To assess the convergence and robustness of the proposed method, a total of $10,000$ initial conditions are randomly generated in a domain 
\textcolor{blue}{$\mathcal{A} := \{ (r_0, v_0, \omega_0,m_0) \mid r_0 \in [1738, 1911.9738]~{\rm km},  
v_0 \in [-83.9779, 83.9779]~\\
{\rm m/s}, \omega_0 \in [0, 9.6638 \times 10^{-4}]~{\rm rad/s},  m_0 \in [240, 600]~{\rm kg}\}$}.
% \begin{equation*}
% \begin{split}
% &r_0 \in [1738, 1911.9738]~({\rm km}), \\
% &v_0 \in [-83.9779, 83.9779]~({\rm m/s}), \\
% &\omega_0 \in [0, 9.6638 \times 10^{-4}]~({\rm rad/s}), \\
% &m_0 \in [240, 600]~({\rm kg}).
% \end{split}
% \end{equation*}
\subsubsection{\textcolor{blue}{Efficacy of Subsection \ref{estimationOF_flighttime} in Section \ref{PIIMTOSLP}}}
\textcolor{blue}{Here we demonstrate the significance of the analytical estimation of the minimum flight time in (\ref{EQ:t_f}). Equation (\ref{EQ:TPBVP_law_ICVN}) is solved for these initial conditions. The corresponding shooting vector is initialized by randomly generating  numbers in the specified domain in (\ref{EQ:solutionspace_ICVN}). On the other hand, we initialize the minimum flight time according to (\ref{EQ:t_f}). Notice that, for a fair comparison, each initial condition among the $10,000$ cases remains the same.}
\begin{table*}[htbp]
  \centering
  \caption{Quantitative comparison of the results by solving $\boldsymbol{\Phi}^f_T(\boldsymbol{z}^{f}_{\text{T,ICVN}})$ in (\ref{EQ:TPBVP_law_ICVN}) using different initial guesses for $t_f$}
  \begin{tabular}{ccc}
  \hline
  Item      & $t_f \in (0, t_\mathrm{max}]$   & $t_f = \hat{t}_f$ \\ 
  \hline
  No. of successful landings & $\textcolor{blue}{5,501}$  &$6,757$  \\ 
  No. of convergent solutions with $t_f < 0$ & $\textcolor{blue}{2,705}$  &$2,677$  \\ 
  No. of convergent solutions with $r(t) < R_0$ & $\textcolor{blue}{80}$  &$75$  \\ 
  Average computational time (s) & $\textcolor{blue}{0.2463}$  &$0.1062$  \\ 
  Average No. of iterations & $\textcolor{blue}{24.87}$  &$21.54$   \\ 
  Average No. of function evaluations & $\textcolor{blue}{154.83}$  &$131.79$  \\ 
  \hline
  \label{Table:t_f_comparison}
  \end{tabular}
  \end{table*}

\textcolor{blue}{
  Table \ref{Table:t_f_comparison} summarizes the numerical results in terms of some important  indicators. To be specific, the successful landing means that not only the shooting function is solved successfully, but also the resulting trajectory is physically feasible, i.e., $t_f >0$ and the radial distance $r(t) \geq R_0$ for $t \in [0,t_f]$. Clearly, providing the analytical estimation of the minimum flight time can improve the success rate for finding the optimal 
  feasible trajectory. Regarding the average computational time, which is calculated as the sum of computational times for successful landings divided by the number of successful landings, initializing $t_f = \hat{t}_f$ reduces the average computational time from $0.2463$ seconds to $0.1062$ seconds. Therefore, we only compare the results obtained by initializing the minimum flight time according to (\ref{EQ:t_f}) hereafter.}
\subsubsection{\textcolor{blue}{Efficacy of Subsections \ref{eliminated_pm0_p0}-\ref{reducing_solution} in Section \ref{PIIMTOSLP}}}
\textcolor{blue}{
  We further solve another two shooting functions in (\ref{EQ:TPBVP_law_ICVN_pm_p0}) and (\ref{EQ:TPBVP_new}) for these initial conditions.
The corresponding shooting vector is initialized according to (\ref{EQ:solutionspace_ICVN_simple}) and (\ref{EQ:shooting_new_back_initial_less}), respectively. 
Table \ref{Table:control_effort} summarizes the numerical results. 
Compared with the indirect methods that propagate the dynamics forward, i.e., solving $\boldsymbol{\Phi}^f_T(\boldsymbol{z}^{f}_{\text{T,ICVN}})$ and  $\boldsymbol{\Phi}^f_T(\boldsymbol{z}^{f}_{\text{T,SICVN}})$, solving $\boldsymbol{\Phi}^b_T(\boldsymbol z^b_{\text{T,SICVN}})$ on the basis of the PIIM presents the highest success rate for finding the optimal feasible trajectory. Although many solutions are convergent, the resulting final times are negative for $\boldsymbol{\Phi}^f_T(\boldsymbol{z}^{f}_{\text{T,ICVN}})$ and  $\boldsymbol{\Phi}^f_T(\boldsymbol{z}^{f}_{\text{T,SICVN}})$, as indicated by $2,677$ and $3,591$ cases. This is because these two shooting functions have a larger solution space than the PIIM-based shooting function. This fact will increase the likelihood of generating a wrong initial guess that leads to the nonlinear function solver converging to a solution with $t_f <0$.}
\begin{table*}[htb]
  \centering
  \caption{Quantitative comparison of the results by solving different shooting functions}
  \begin{tabular}{ccc}
  \hline
  Item          & $\boldsymbol{\Phi}^f_T(\boldsymbol{z}^{f}_{\text{T,SICVN}})$ \textcolor{black}{in (\ref{EQ:TPBVP_law_ICVN_pm_p0})}  & $\boldsymbol{\Phi}^b_T(\boldsymbol z^b_{\text{T,SICVN}})$ \textcolor{black}{in (\ref{EQ:TPBVP_new})}\\ 
  \hline
  No. of successful landings   &$5,731$  &$9,887$   \\ 
  No. of convergent solutions with $t_f < 0$   &$3,591$  &$42$  \\ 
  No. of convergent solutions with $r(t) < R_0$   &$92$  &$48$    \\ 
  Average computational time (s)   &$0.0831$  &$0.0613$   \\ 
  Average No. of iterations   &$30.79$  &$13.49$   \\ 
  Average No. of function evaluations   &$128.46$  &$68.58$ \\ 
  \hline
  \label{Table:control_effort}
  \end{tabular}
  \end{table*}
  \begin{table*}[hbt!]
    \centering
    \caption{Quantitative comparison of the results after using the remedy strategy}
    \begin{tabular}{cccc}
    \hline
    Item      & $\boldsymbol{\Phi}^f_T(\boldsymbol{z}^{f}_{\text{T,ICVN}})$ \textcolor{black}{in (\ref{EQ:TPBVP_law_ICVN})}     & $\boldsymbol{\Phi}^f_T(\boldsymbol{z}^{f}_{\text{T,SICVN}})$ \textcolor{black}{in (\ref{EQ:TPBVP_law_ICVN_pm_p0})}  & $\boldsymbol{\Phi}^b_T(\boldsymbol z^b_{\text{T,SICVN}})$ \textcolor{black}{in (\ref{EQ:TPBVP_new})}\\  
    \hline
    No. of successful landings & $9,928$  &$9,930$  &$9,952$   \\ 
    No. of convergent solutions with $t_f < 0$ & $0$  &$0$  &$0$  \\ 
    No. of convergent solutions with $r(t) < R_0$ & $48$  &$48$  &$48$    \\ 
    Average computational time (s) & $0.0903$  &$0.0955$  &$0.0557$   \\ 
    Average No. of iterations & $16.53$  &$34.87$  &$11.21$   \\ 
    Average No. of function evaluations & $102.69$  &$148.00$  &$58.71$ \\ 
    \textcolor{black}{Success rate ($\%$)}  & $99.76$  & $99.78$  & $100$ \\
    \hline
    \label{Table:control_effort_remedy}
    \end{tabular}
    \end{table*}

Furthermore, these three shooting functions result in $75$, $92$, and $48$ cases of convergent but infeasible solutions with $r(t) < R_0$ (indicating that the lunar lander is flying underneath the lunar surface), respectively. For the average computational time, solving $\boldsymbol{\Phi}^b_T(\boldsymbol z^b_{\text{T,SICVN}})$ in (\ref{EQ:TPBVP_new}) requires the shortest time compared to other two. Meanwhile, the numbers of iterations and function evaluations are also the lowest among the three.
% \begin{table*}[htbp]
% \centering
% \caption{Quantitative comparison of the results by solving different shooting functions}
% \begin{tabular}{cccc}
% \hline
% Item      & $\boldsymbol{\Phi}^f_T(\boldsymbol{z}^{f}_{\text{T,ICVN}})$ \textcolor{black}{in (\ref{EQ:TPBVP_law_ICVN})}     & $\boldsymbol{\Phi}^f_T(\boldsymbol{z}^{f}_{\text{T,SICVN}})$ \textcolor{black}{in (\ref{EQ:TPBVP_law_ICVN_pm_p0})}  & $\boldsymbol{\Phi}^b_T(\boldsymbol z^b_{\text{T,SICVN}})$ \textcolor{black}{in (\ref{EQ:TPBVP_new})}\\ 
% \hline
% No. of successful landings & $6,757$  &$5,731$  &$9,887$   \\ 
% No. of convergent solutions with $t_f < 0$ & $2,677$  &$3,591$  &$42$  \\ 
% No. of convergent solutions with $r(t) < R_0$ & $75$  &$92$  &$48$    \\ 
% Average computational time (s) & $0.1062$  &$0.0831$  &$0.0613$   \\ 
% Average No. of iterations & $21.54$  &$30.79$  &$13.49$   \\ 
% Average No. of function evaluations & $131.79$  &$128.46$  &$68.58$ \\ 
% \hline
% \label{Table:control_effort}
% \end{tabular}
% \end{table*}
\subsubsection{\textcolor{blue}{A remedy strategy ensuring $t_f > 0$}}
Notice that the initial guess of the minimum flight time remains the same for all three shooting functions, so it is reasonable to conclude that providing an inappropriate initial guess of the co-state vector increases the likelihood of converging to an infeasible solution with a negative final time, as indicated by the numbers of convergent solutions with $t_f < 0$ in Table \ref{Table:control_effort}. To resolve this issue, a simple remedy strategy is proposed here. 

\textcolor{black}{
The variable $t_f$ in the shooting function is first replaced by a new variable $\xi$. Then the integration interval is changed from $[0,t_f]$ to $[0, \exp(\xi)]$. This ensures that the iteration of the final time will remain positive. It should be
noted that the initial guess of $\xi$ is given by $\log(\hat{t}_f)$, where $\hat{t}_f$
is calculated according to (\ref{EQ:t_f}).}

\begin{figure*}[htb!] 
  \begin{subfigure}{0.48\textwidth}
  \includegraphics[width=\textwidth]{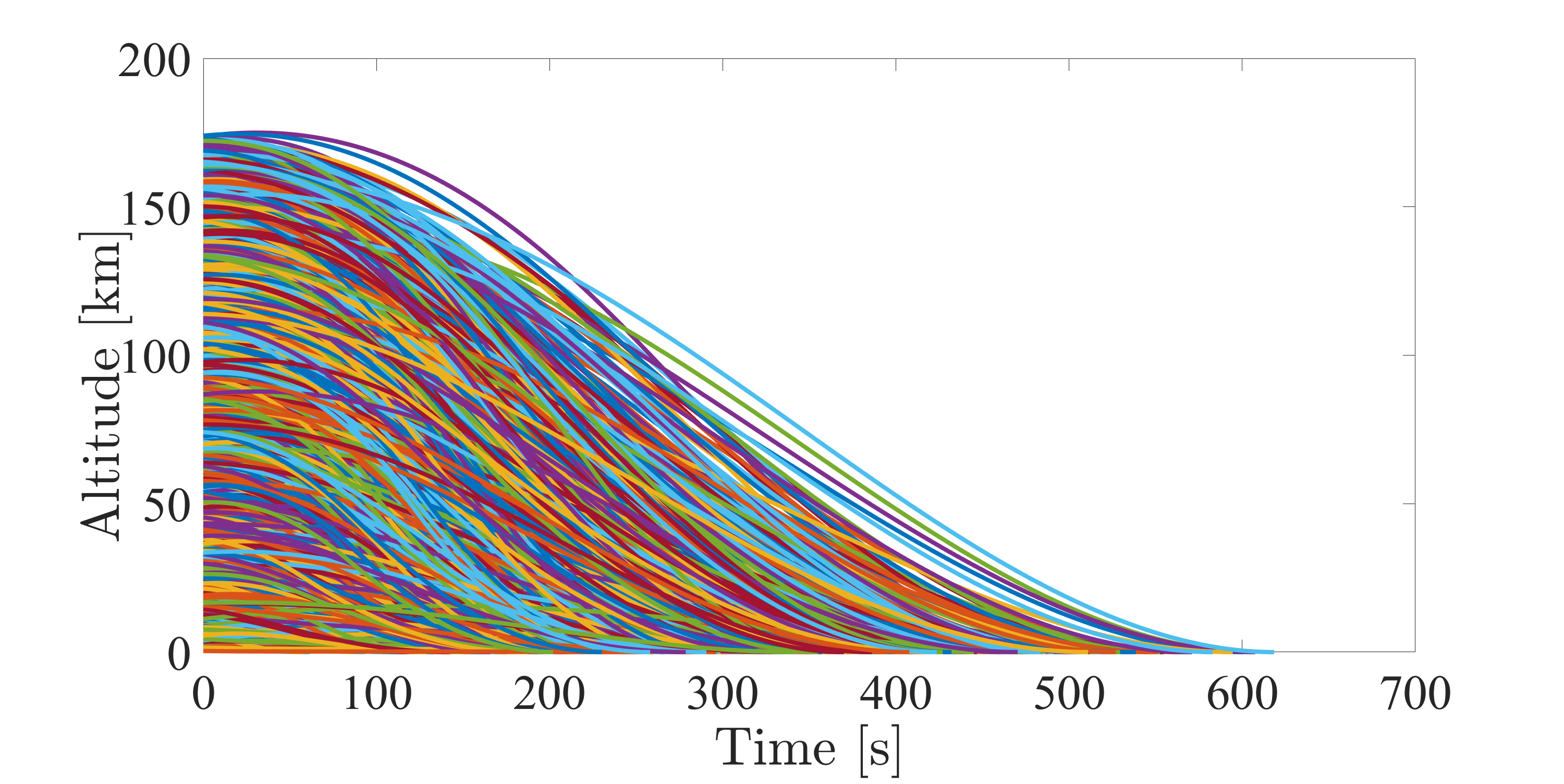}
  \caption{Altitude profiles.}
  \label{Fig:cooperative_control_1}
  \end{subfigure}
  \begin{subfigure}{0.48\textwidth}
  \includegraphics[width=\textwidth]{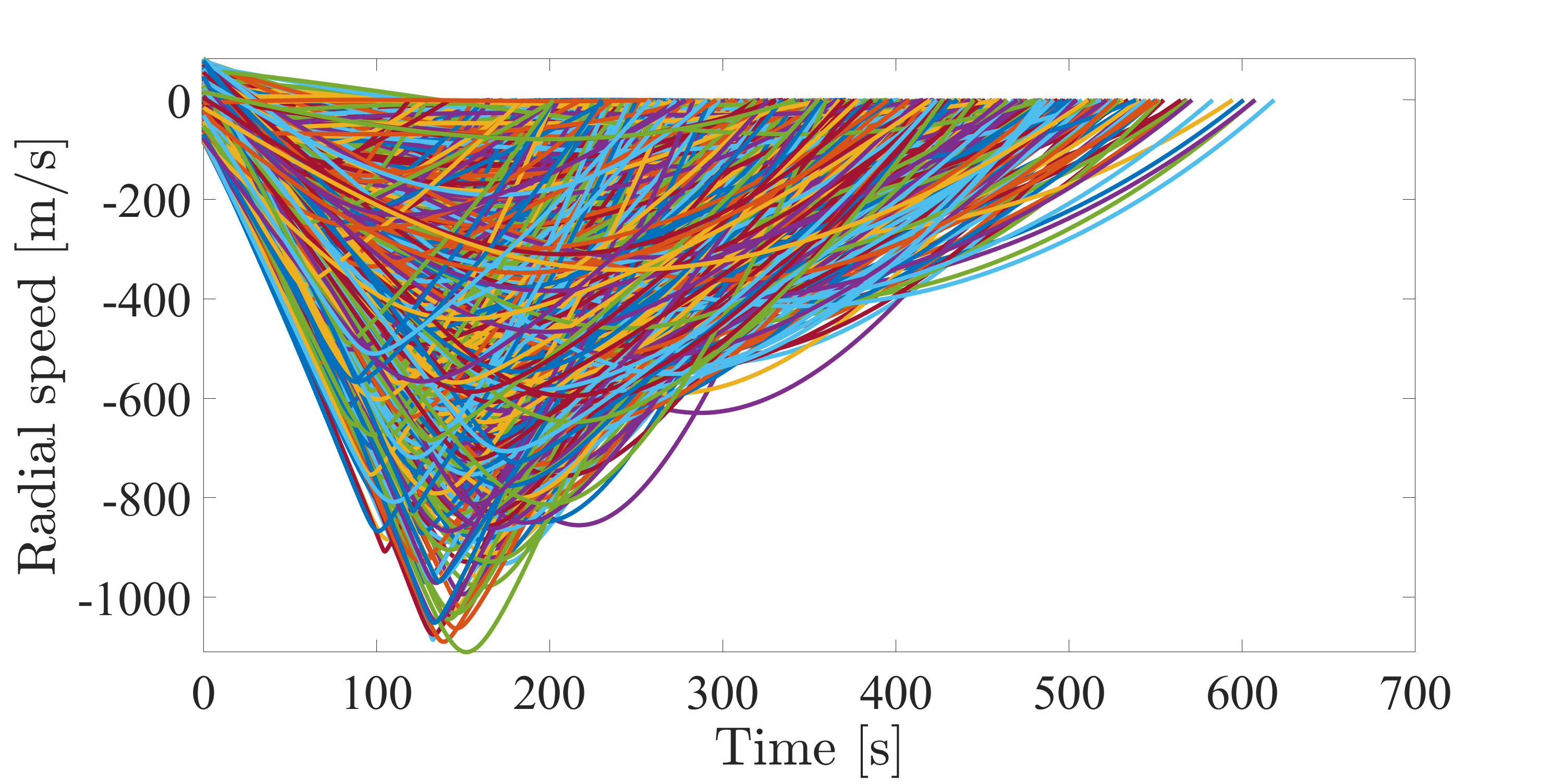}
  \caption{Radial speed profiles.}
  \label{Fig:cooperative_control_2}
  \end{subfigure}
  \begin{subfigure}{0.48\textwidth}
  \includegraphics[width=\textwidth]{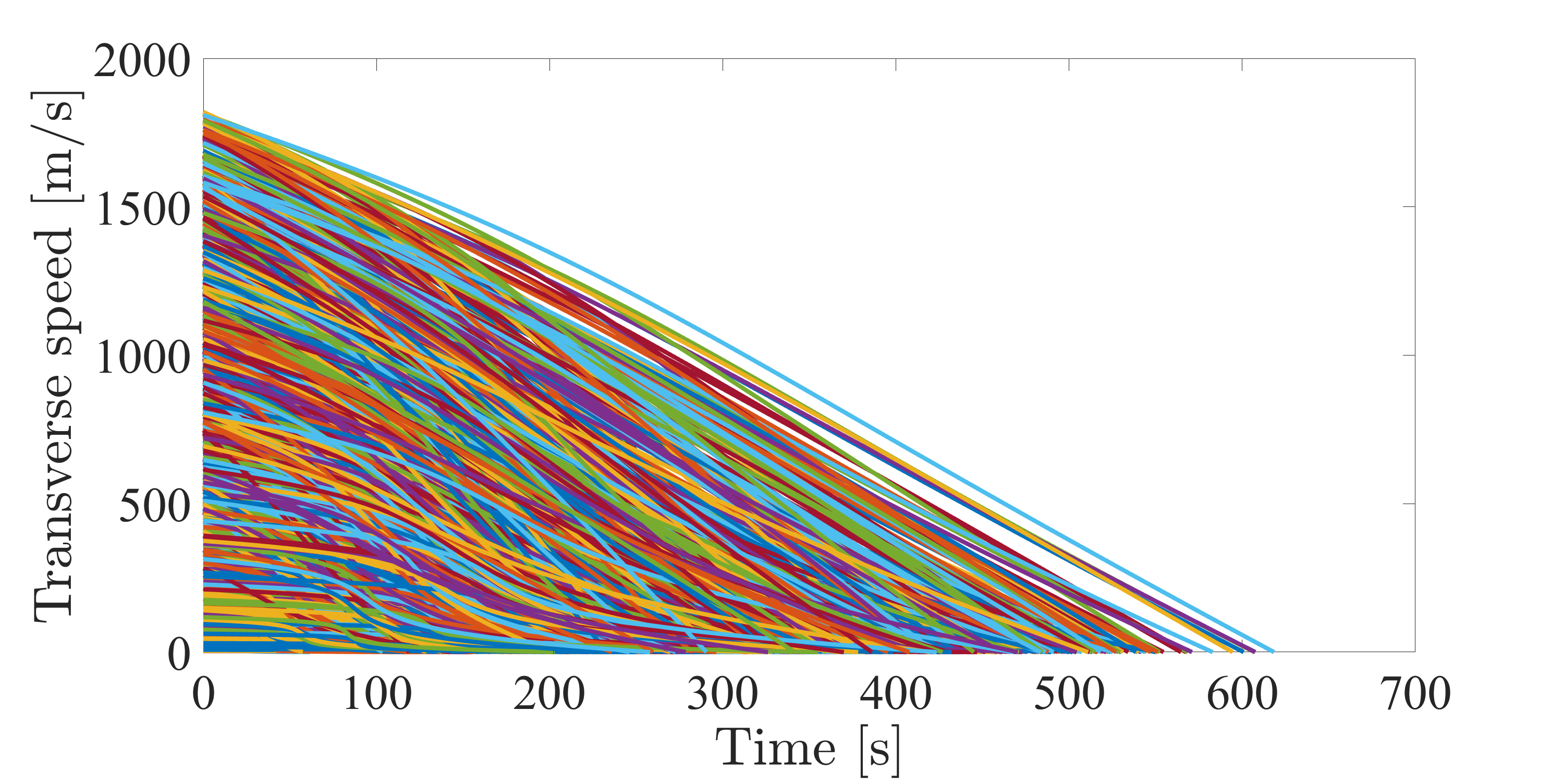}
  \caption{Transverse speed profiles.}
  \label{Fig:cooperative_control_3}
  \end{subfigure}
  \begin{subfigure}{0.48\textwidth}
  \includegraphics[width=\textwidth]{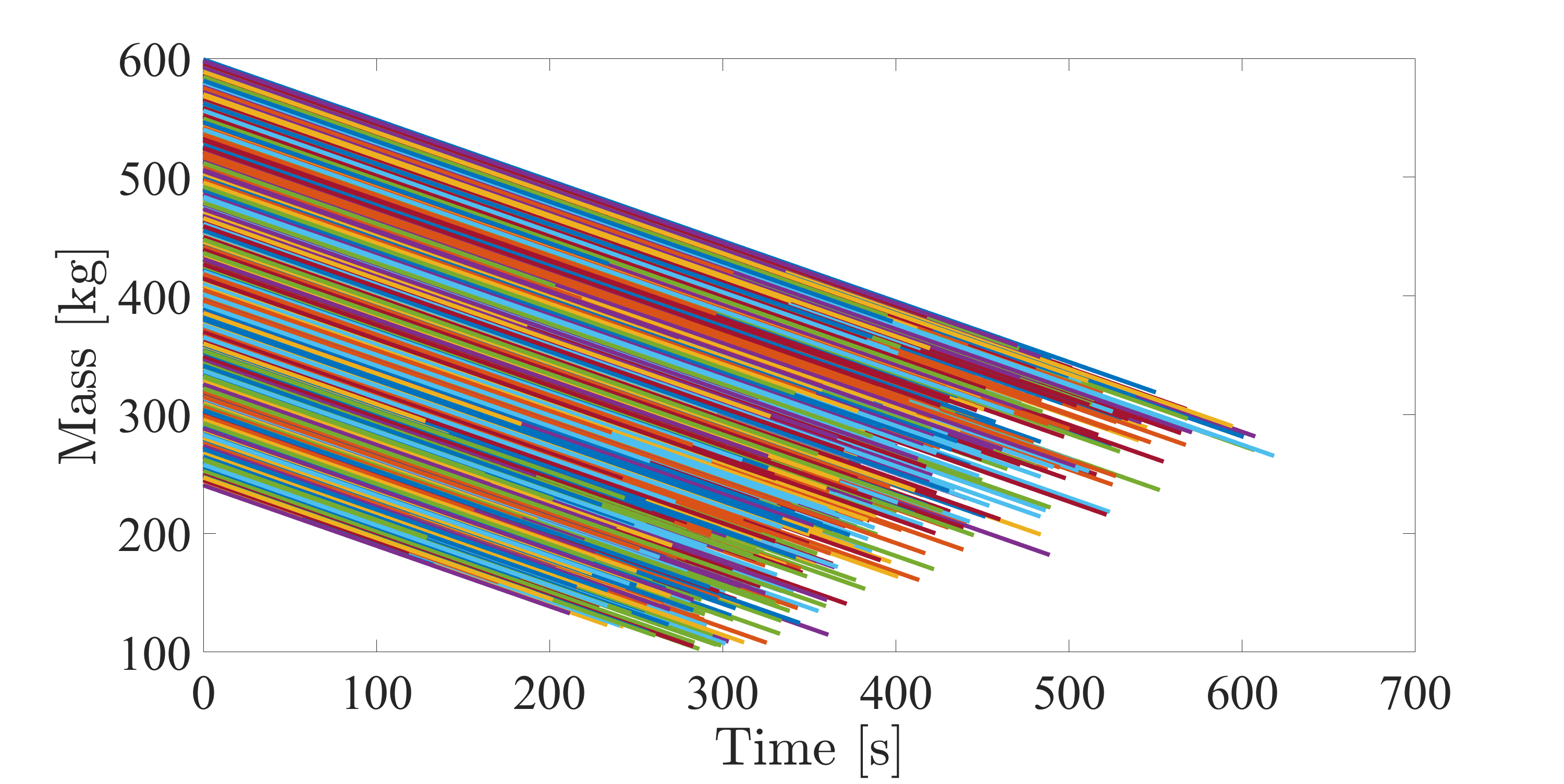}
  \caption{Mass profiles.}
  \label{Fig:cooperative_control_4}
  \end{subfigure}
  \caption{Profiles of the altitude, radial speed, transverse speed, and mass.}
  \label{Fig:cooperative_profile}
  \end{figure*}
With the remedy strategy, we run another test for the same $10,000$ cases. The numerical results are displayed in Table \ref{Table:control_effort_remedy}, from which we can see that, the numbers of successful landings are improved for the three shooting functions. The reported numbers of convergent solutions with $t_f < 0$ are all reduced to zero as expected. Meanwhile, they share the same number of convergent solutions with $r(t) < R_0$, which is $48$ out of $10,000$ cases. 
Since the initial condition is randomly generated in $\mathcal{A}$, it is possible that the initial condition does not have a solution inherently. For the $48$ cases, we find that $r(t) < R_0$ holds for all obtained solutions despite trying a large number of initial guesses, so it is reasonable to conclude that these cases do not have feasible solutions.
In such case, we can see that solving $\boldsymbol{\Phi}^b_T(\boldsymbol z^b_{\text{T,SICVN}})$ in (\ref{EQ:TPBVP_new}) exhibits a success rate of $100\%$. In contrast, solving $\boldsymbol{\Phi}^f_T(\boldsymbol{z}^{f}_{\text{T,ICVN}})$ and $\boldsymbol{\Phi}^f_T(\boldsymbol{z}^{f}_{\text{T,SICVN}})$ results in a success rate of $99.76\%$ and $99.78\%$, respectively.
Moreover, after implementing the remedy strategy, the performance of solving $\boldsymbol{\Phi}^b_T(\boldsymbol z^b_{\text{T,SICVN}})$ is further improved and outperforms the other two shooting functions, as indicated by the average computational time, number of iterations and function evaluations in Table \ref{Table:control_effort_remedy}.  
Through solving $\boldsymbol{\Phi}^b_T(\boldsymbol z^b_{\text{T,SICVN}})$, the average computational time is decreased by $38.32\%$ and $41.68\%$ compared with $\boldsymbol{\Phi}^f_T(\boldsymbol{z}^{f}_{\text{T,ICVN}})$ and $\boldsymbol{\Phi}^f_T(\boldsymbol{z}^{f}_{\text{T,SICVN}})$, respectively. Surprisingly, although there are two less unknown variables in $\boldsymbol{\Phi}^f_T(\boldsymbol{z}^{f}_{\text{T,SICVN}})$ than that in $\boldsymbol{\Phi}^f_T(\boldsymbol{z}^{f}_{\text{T,ICVN}})$, the solving performance is not necessarily improved in terms of the success rate and the average computational time. More importantly, compared with a computational time of nearly $1$ second in \cite{pengkun2019} where an implicit shooting method was proposed, it only takes $0.0557$ seconds to find the optimal solution by solving $\boldsymbol{\Phi}^b_T(\boldsymbol z^b_{\text{T,SICVN}})$. \textcolor{blue}{Thanks to the PIIM, the computational time and success rate of solving the TOSLP are improved by $77.39\%$ and $44.72\%$, respectively over the conventional indirect method in Section \ref{SE:indirect_TOSLP}, as shown by the results in Tables \ref{Table:t_f_comparison} and \ref{Table:control_effort_remedy}.}

\begin{figure}[htb]
\begin{center}
\includegraphics[scale=0.2]{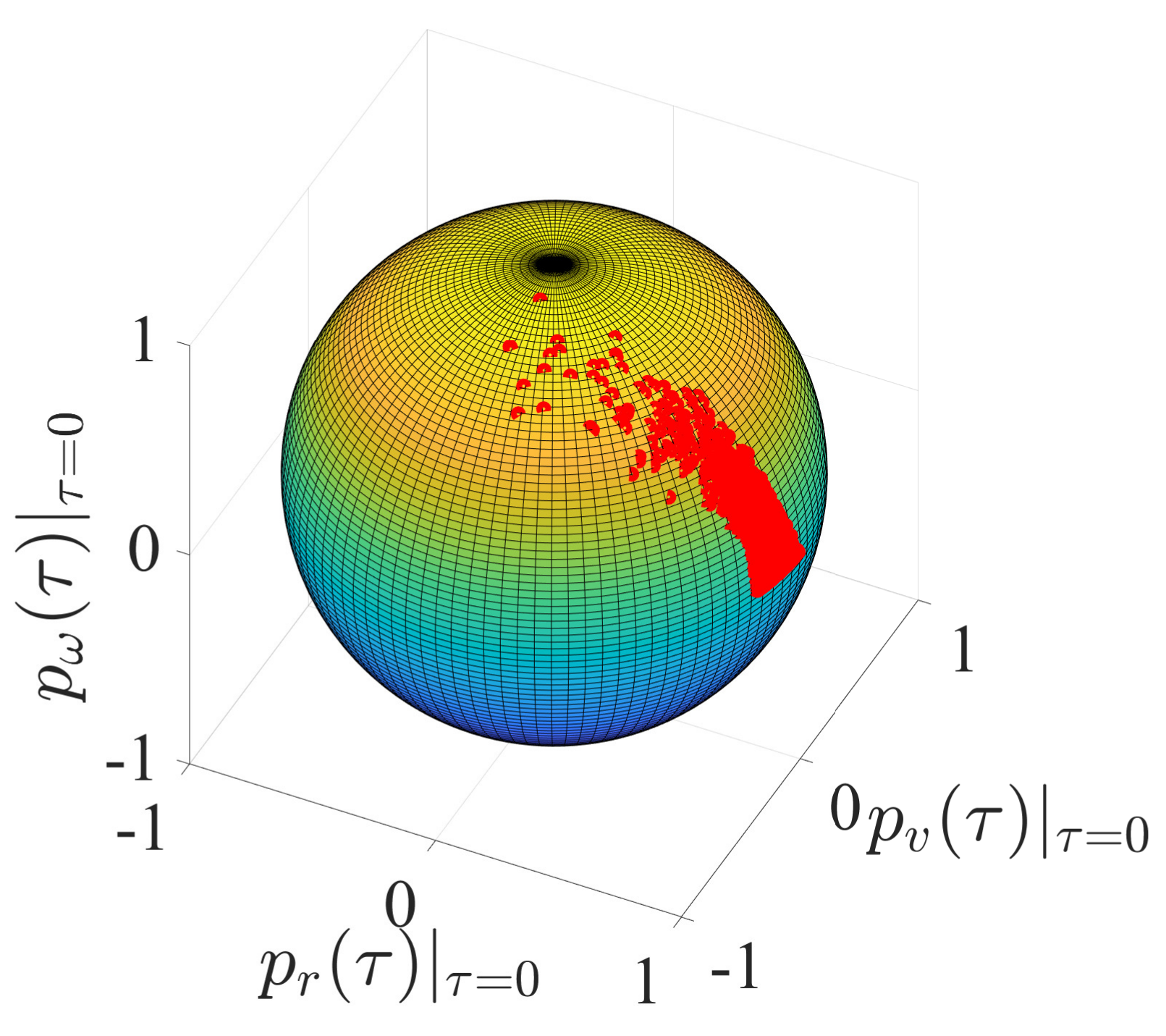}
\caption{Distribution of  $p_r(\tau)|_{\tau=0}$, $p_v(\tau)|_{\tau=0}$, and $p_\omega(\tau)|_{\tau=0}$.}\label{Fig:estimate_costatesphere}
\end{center}
\end{figure}

\begin{figure}[htb]
  \begin{center}
  \includegraphics[scale=0.2]{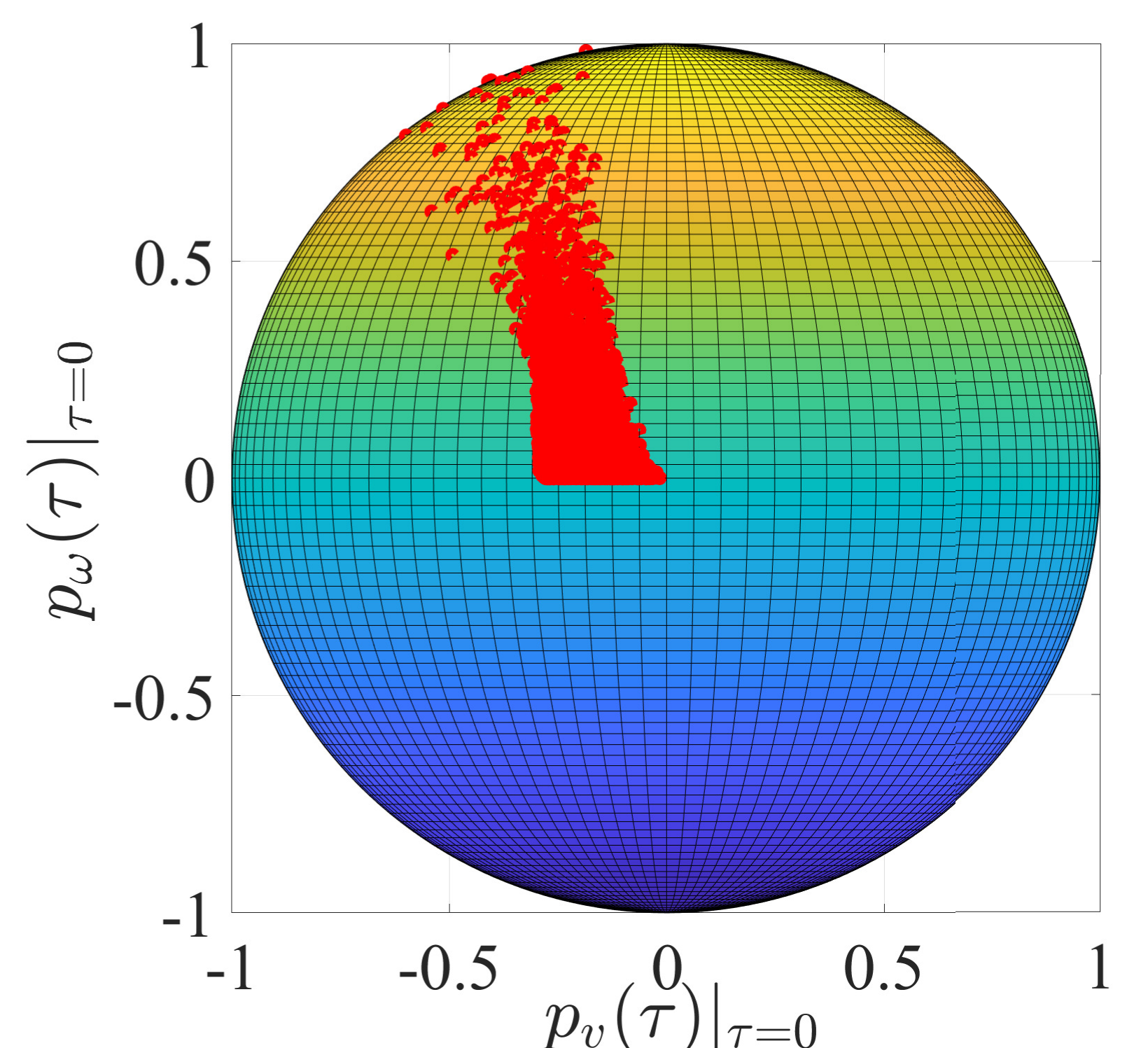}
  \caption{Distribution of $p_v(\tau)|_{\tau=0}$ and $p_\omega(\tau)|_{\tau=0}$.}\label{Fig:estimate_costate_planar}
  \end{center}
  \end{figure}

Out of the $9,952$ successful landing cases, a total of $1,000$ landings are displayed in Fig.~\ref{Fig:cooperative_profile} in terms of altitude, radial speed, transverse speed, and mass.
Although the initial condition is randomly generated, solving $\boldsymbol{\Phi}^b_T(\boldsymbol z^b_{\text{T,SICVN}})$ can always find the optimal solution in a fast and reliable manner. 
To further validate the initial guess of $p_r(\tau)|_{\tau=0}$, $p_v(\tau)|_{\tau=0}$, and $p_\omega(\tau)|_{\tau=0}$ defined in (\ref{EQ:shooting_new_back_initial_less}), 
the convergent solutions to $p_r(\tau)|_{\tau=0}$, $p_v(\tau)|_{\tau=0}$, and $p_\omega(\tau)|_{\tau=0}$ for the $9,952$ successful landing cases are shown in Figs.~\ref{Fig:estimate_costatesphere} and \ref{Fig:estimate_costate_planar}. 
It can be observed that  $p_r(\tau)|_{\tau=0} > 0$ holds true as expected by (\ref{EQ:end_pr_positive}). 
Meanwhile, (\ref{EQ:end_pv_pw}) can also be validated, as shown by Fig.~\ref{Fig:estimate_costate_planar}, in which $p_r(\tau)|_{\tau=0}$, $p_v(\tau)|_{\tau=0}$, and $p_\omega(\tau)|_{\tau=0}$ are constrained on a unit 3-D octant sphere. 
\subsection{Simulations on the FOSLP}\label{SimB}
In the \textcolor{blue}{preceding} subsection, we demonstrated that fast and robust convergence can be achieved for the TOSLP 
by solving $\boldsymbol{\Phi}^b_T(\boldsymbol z^b_{\text{T,SICVN}})$. Now, we will show that the developments of the \textcolor{black}{PIIM} for solving the TOSLP, can also be applied to facilitate convergence \textcolor{black}{of} the FOSLP. 

\begin{figure}[hbt!]      
  \begin{center}
  \includegraphics[scale=0.2]{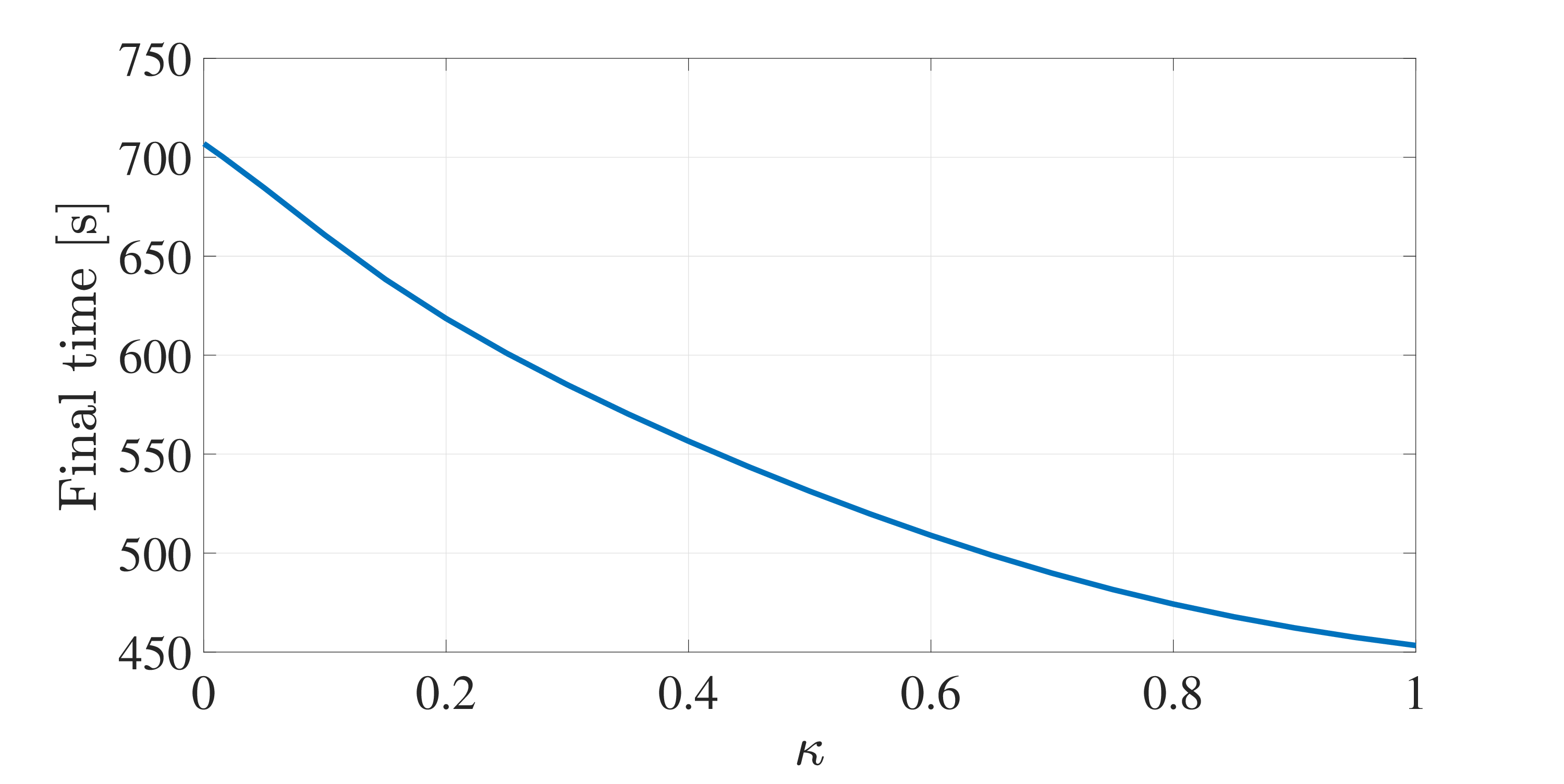}
  \caption{Homotopy path for $t_f(\kappa)$ with $\delta = 0.1$.}\label{Fig:kappa_tf}
  \end{center}
  \end{figure}

To illustrate the homotopy process, the following initial condition is considered: $r_0=1,902.1754$ km, $v_0 = 23.1290$ m/s, $\omega_0 = 2.3261 \times 10^{-4}$ rad/s, and $m_0=483.4040$ kg.
By solving (\ref{EQ:TPBVP_new}) for the TOSLP, $p_0$ is found to be $0.5693$, which will be used and kept unchanged in (\ref{EQ:performance_homotopy}) throughout the homotopy process. Fig.~\ref{Fig:kappa_tf} plots the homotopy path of $t_f(\kappa)$, which shows the profile of the convergent solution to the optimal final time $t_f$ as the homotopy parameter $\kappa$ is decreased from 1 to 0. 
During the solution process, we observe that very small values of $\kappa$ may impair convergence. Thus, the homotopy process for $\kappa$ is terminated once $\kappa$ is decreased to \textcolor{black}{a number less than} $2^{-4}$. It is worth mentioning that $\delta$ is kept as $0.1$ before the homotopy process for $\kappa$ is completed. From Fig.~\ref{Fig:kappa_tf}, the homotopy path turns out to be continuous and unidirectional, which makes the homotopy process straightforward \cite{pan2016double}. 

\begin{figure}[hbt!]
\begin{center}
\includegraphics[scale=0.2]{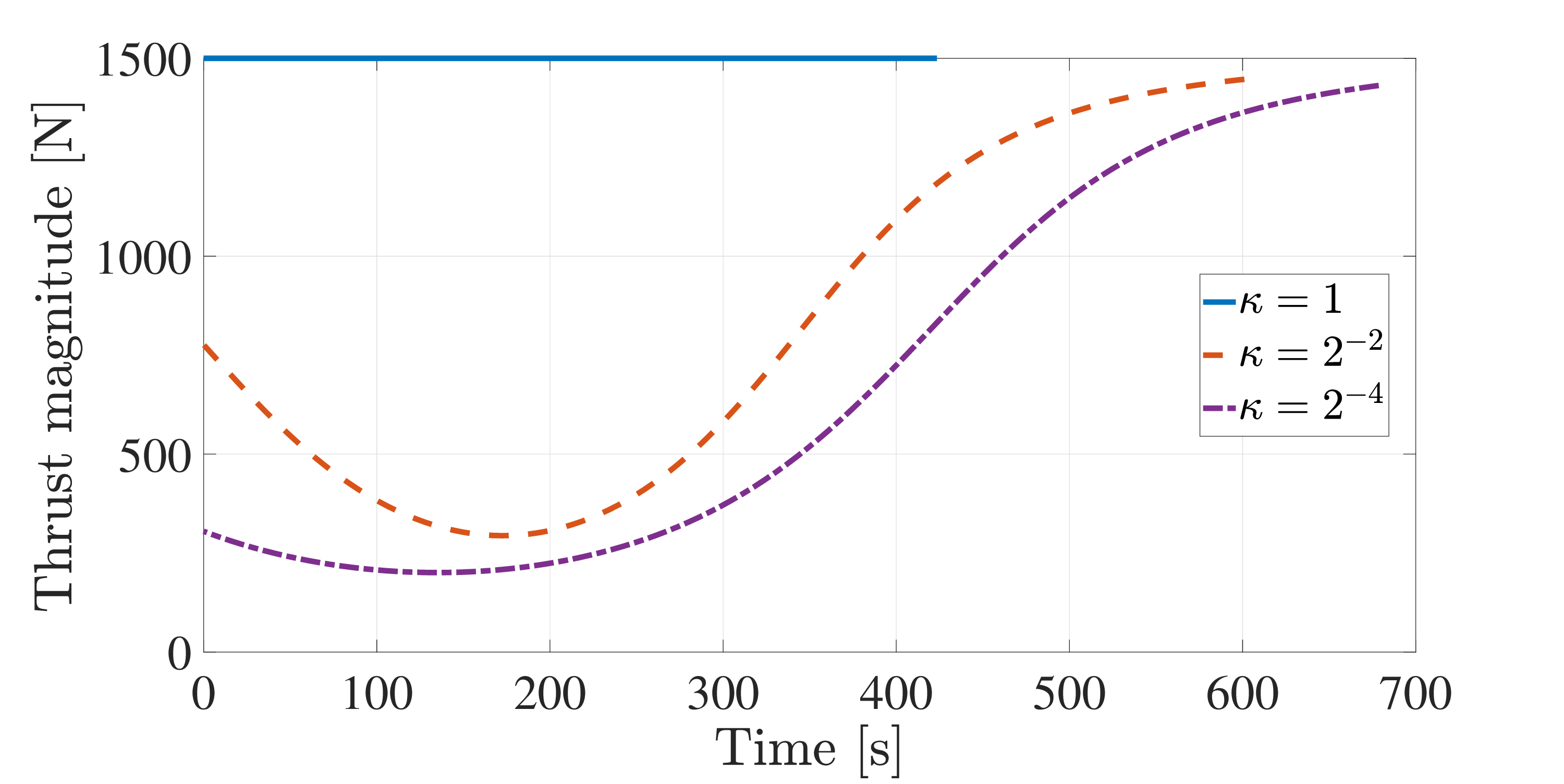}
\caption{Thrust magnitude profiles with different $\kappa$.}\label{Fig:kappa_u}
\end{center}
\end{figure}

\begin{figure}[hbt!]
\begin{center}
\includegraphics[scale=0.2]{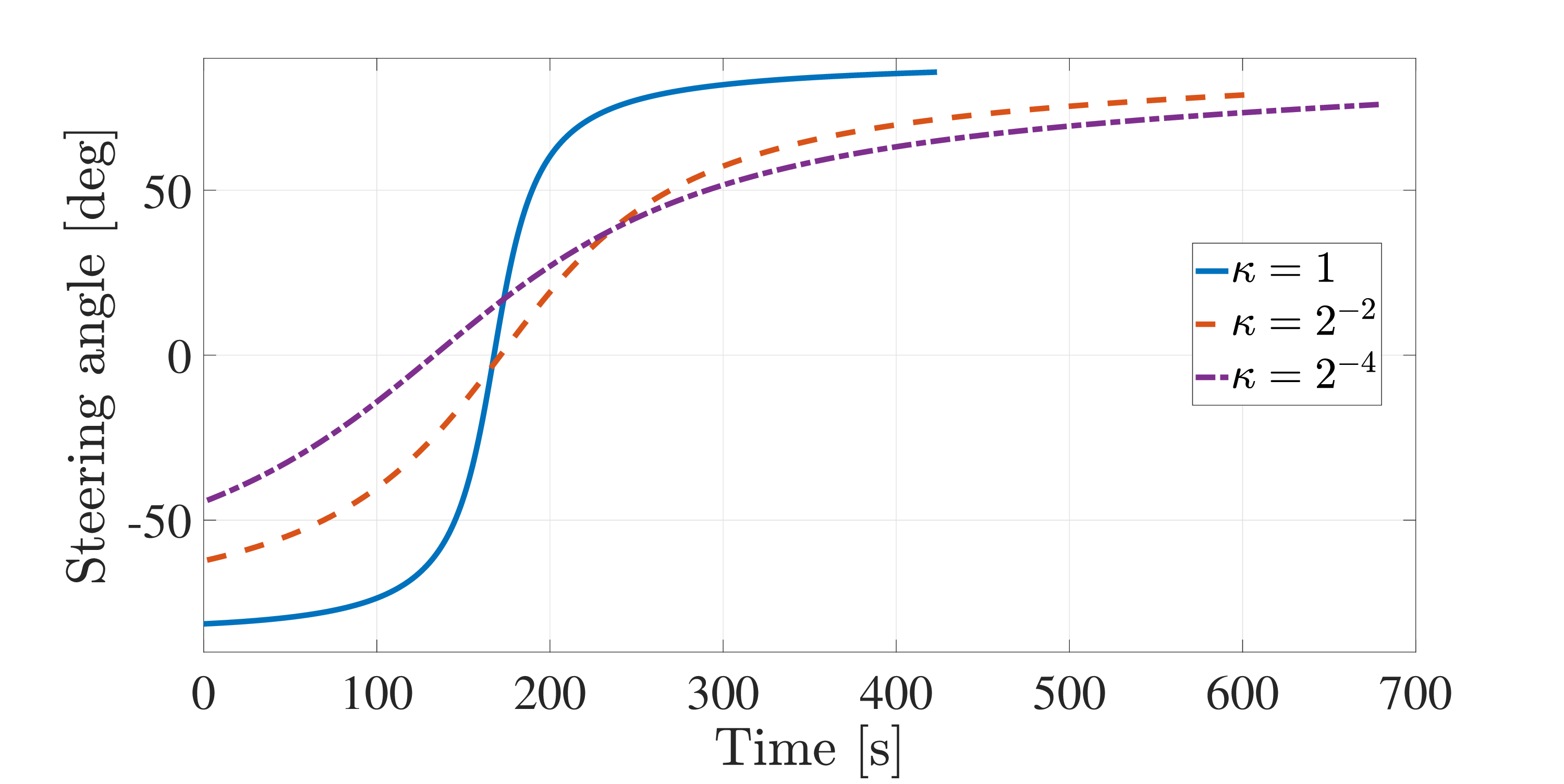}
\caption{Thrust steering angle profiles with different $\kappa$.}\label{Fig:kappa_angle}
\end{center}
\end{figure}

\begin{figure}[hbt!]
\begin{center}
\includegraphics[scale=0.2]{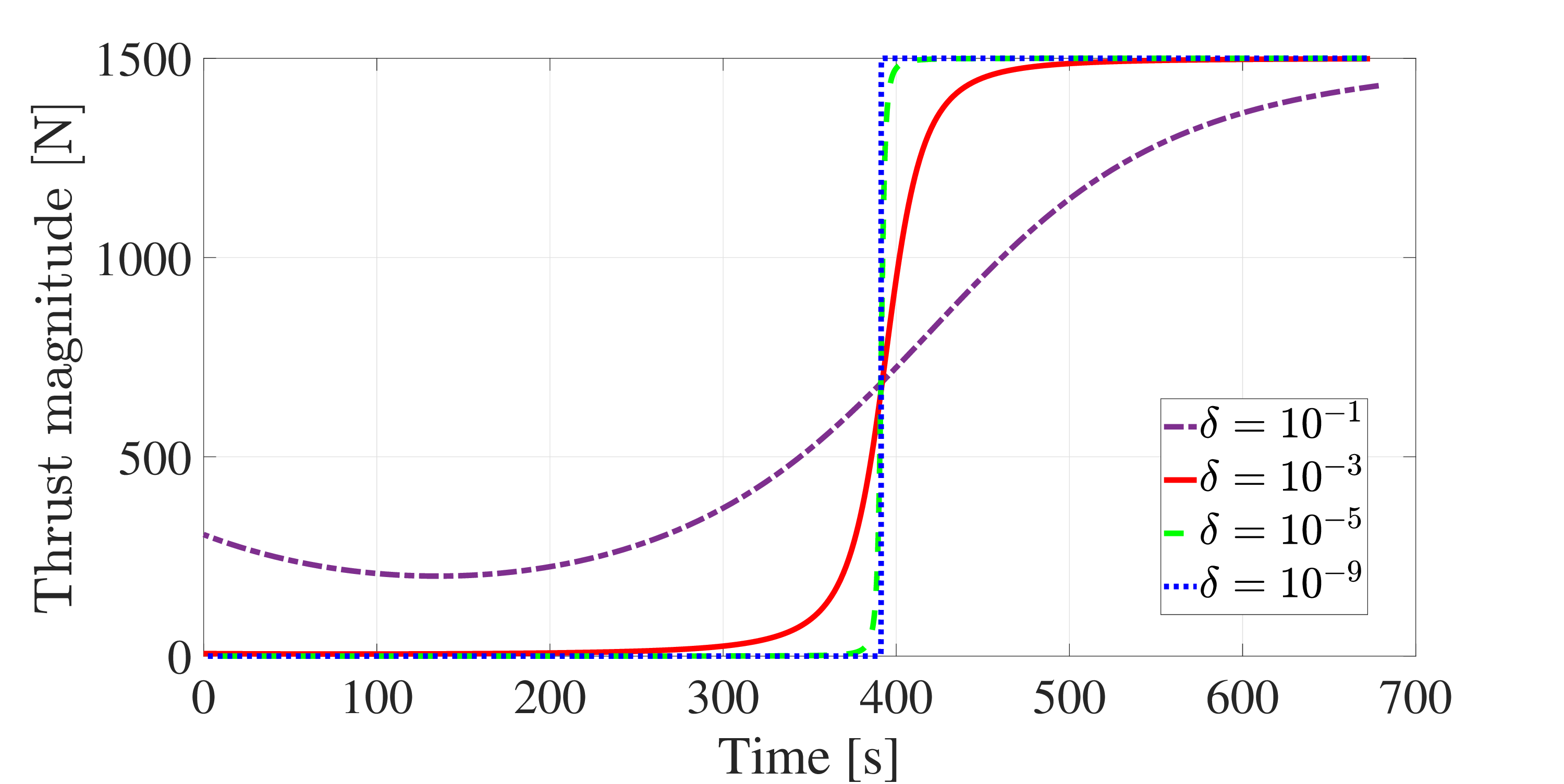}
\caption{Thrust magnitude profiles with different $\delta$ and $\kappa = 2^{-4}$.}\label{Fig:alpha_u}
\end{center}
\end{figure}

\begin{figure}[hbt!]
\begin{center}
\includegraphics[scale=0.2]{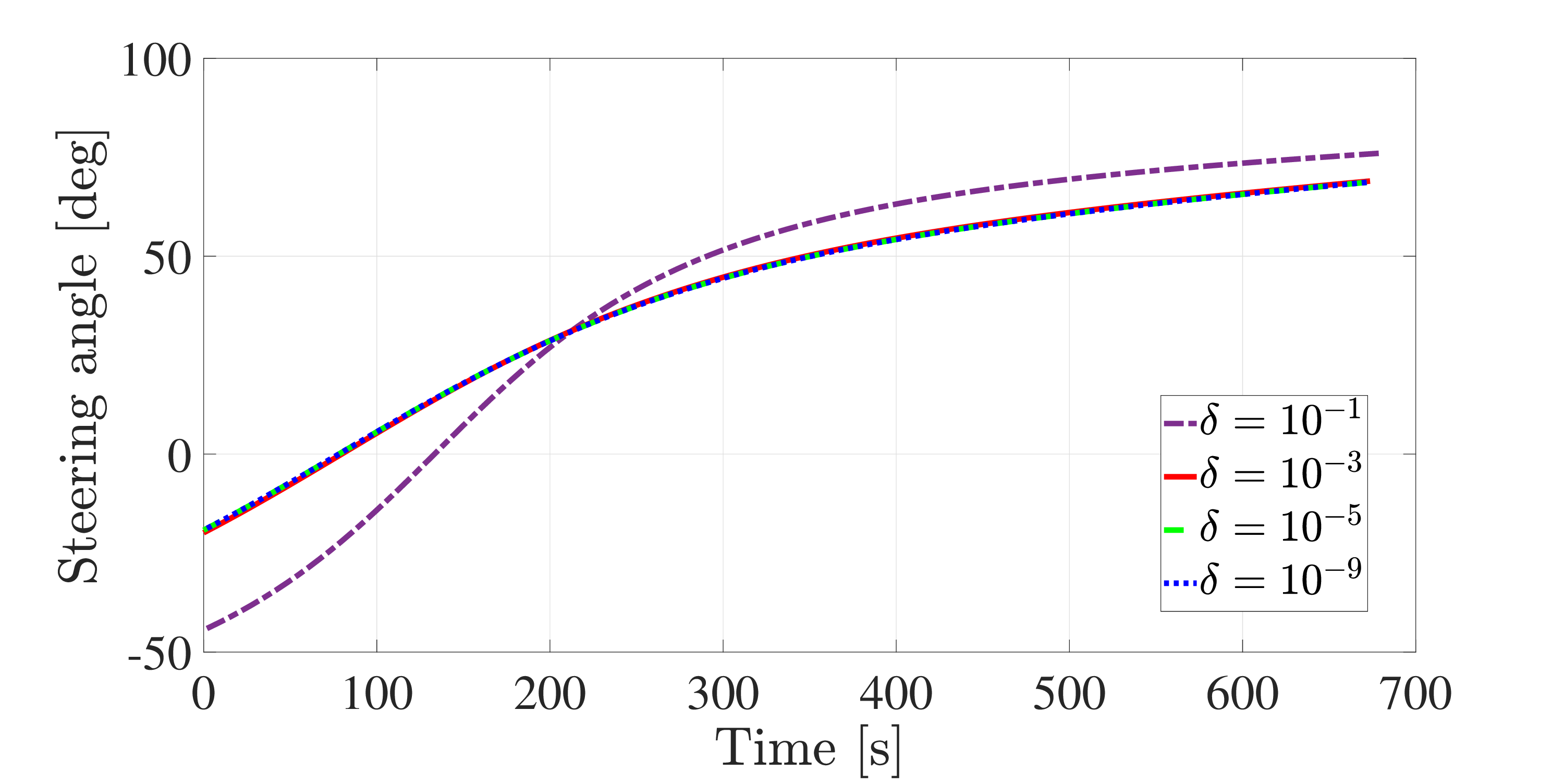}
\caption{Thrust steering angle profiles with different $\delta$ and $\kappa = 2^{-4}$.}\label{Fig:alpha_angle}
\end{center}
\end{figure}

\begin{figure}[hbt!]
\begin{center}
\includegraphics[scale=0.2]{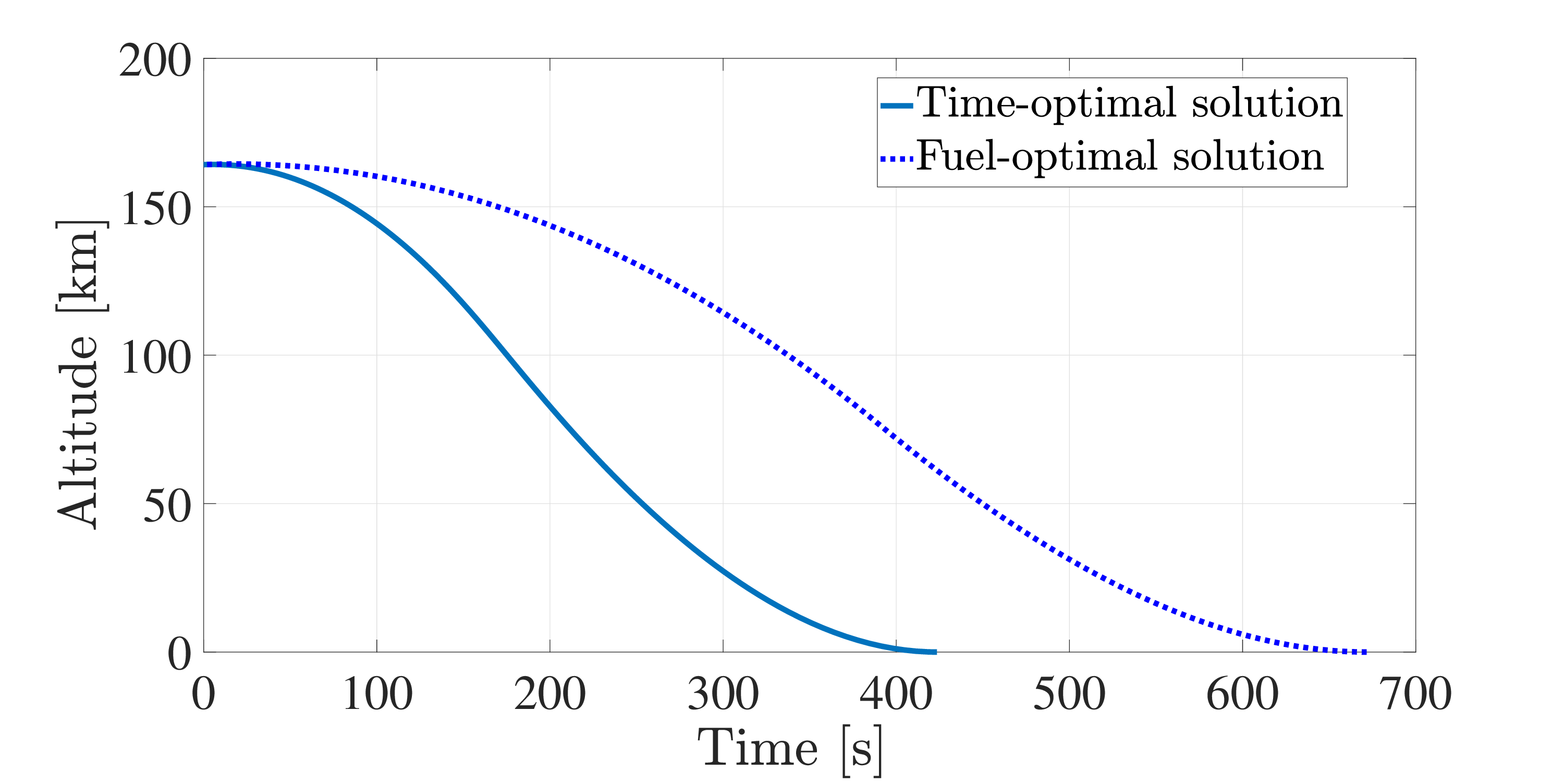}
\caption{Altitude profiles for the TOSLP and the FOSLP.}\label{Fig:altitude_angle}
\end{center}
\end{figure}

\begin{figure}[hbt!]
\begin{center}
\includegraphics[scale=0.2]{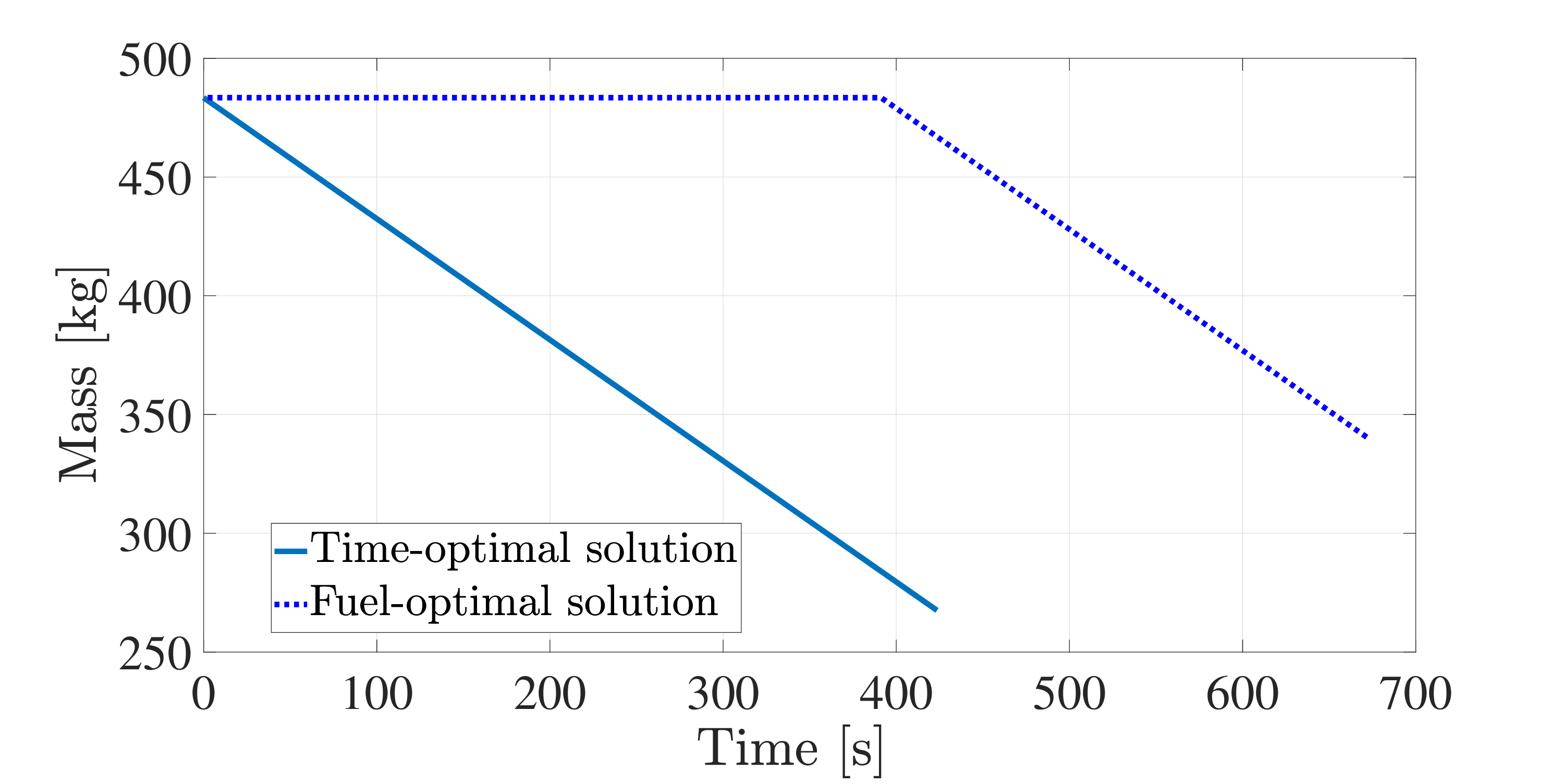}
\caption{Mass profiles for the TOSLP and the FOSLP.}\label{Fig:mass_angle}
\end{center}
\end{figure}

The optimal thrust magnitude and thrust steering angle profiles are displayed in Figs.~\ref{Fig:kappa_u} and \ref{Fig:kappa_angle}, respectively. The solution with $\kappa=1$ represents the optimal solution to the TOSLP. Once $\kappa$ is reduced \textcolor{black}{below} $2^{-4}$, the homotopy process for $\delta$ commences. Figs.~\ref{Fig:alpha_u} and \ref{Fig:alpha_angle} illustrate the optimal thrust magnitude and thrust steering angle profiles as the smoothing constant $\delta$ decreases from $0.1$ to $10^{-9}$. It can be observed that while the intermediate thrust magnitude profiles are continuous, the final thrust magnitude profile with $\delta = 10^{-9}$ exhibits a bang-bang solution with a single switch.
To provide a comprehensive analysis, the optimal altitude and mass profiles related to the TOSLP and the FOSLP are shown in Figs.~\ref{Fig:altitude_angle} and \ref{Fig:mass_angle}, respectively.

To be specific, the minimum flight time is determined to be $423.483$ seconds with a fuel consumption of $215.842$ kg, and the thrust magnitude is kept at the maximum during landing. On the other hand, the final time for the FOSLP is $671.638$ seconds with a fuel consumption of $142.905$ kg, and the thrust magnitude is a typical ``off-on'' form. For comparison, we utilize the conventional indirect method, which resolves $\boldsymbol{\Phi}^f_F(\boldsymbol{z}^{f}_{\text{F,ICVN}}(\delta))$, to determine the optimal fuel consumption. With the same smoothing constant $\delta = 10^{-9}$, the fuel consumption derived from $\boldsymbol{\Phi}^f_F(\boldsymbol{z}^{f}_{\text{F,ICVN}}(\delta))$ is found to be $142.900$ kg, which indicates that the homotopy parameter $\kappa = 2^{-4}$ incurs only a negligible penalty of $0.005$ kg.

To illustrate the robustness, we compare the homotopy approach established in Fig.~\ref{Fig:homo_flowchart} with the conventional method by directly solving (\ref{EQ:TPBVP_law_ICVN_fuel}), in which the initial guess of the \textcolor{black}{shooting vector} is generated according to (\ref{EQ:solutionspace_ICVN_fuel}) and the final time $t_f$ is initialized according to (\ref{EQ:t_f}). The $9,952$ successful landing cases from the preceding subsection are tested. According to \cite{yang2017rapid}, if the TOSLP has a feasible optimal solution, the FOSLP is also feasible under the same initial condition. 

\textcolor{blue}{To further demonstrate the advantages of propagating the dynamics backward, we consider the same homotopy method as in (\ref{EQ:TPBVP_new_homo}) to solve the FOSLP except that the dynamics are propagated forward.
 Notice that the solution to (\ref{EQ:TPBVP_law_ICVN_pm_p0}) can also be used to solve the FOSLP via a homotopy approach. In such case, the dynamics of both states and co-states are propagated forward, and the shooting function in (\ref{EQ:TPBVP_new_homo}) becomes
\begin{equation}
  \boldsymbol{\Phi}^f_h(\boldsymbol z^f_{\text{h,SICVN}}(\kappa,\delta)) = [r(\kappa,\delta)-R_0;v(\kappa,\delta); \omega(\kappa,\delta);p_m(\kappa,\delta);\mathscr{H}_h(\kappa,\delta)]|_{t_f}=\boldsymbol{0},
  \label{EQ:TPBVP_new_homo_forwaed}
  \end{equation}
  where 
$\boldsymbol z^f_{\text{h,SICVN}}(\kappa,\delta) = [p_{r_0}(\kappa,\delta),p_{v_0}(\kappa,\delta), p_{\omega_0}(\kappa,\delta), 
p_{m_0}(\kappa,\delta), t_f(\kappa,\delta)]^T    
$. Once the solution to (\ref{EQ:TPBVP_law_ICVN_pm_p0}) is obtained, the similar process in Fig.~\ref{Fig:homo_flowchart} is implemented. The solutions to the  $9,930$ successful landings obtained by solving (\ref{EQ:TPBVP_law_ICVN_pm_p0}) are employed to solve (\ref{EQ:TPBVP_new_homo_forwaed}). It is worth noting that the value for $p_{m_0}$ is obtained according to (\ref{EQ:deltapm}).}

Table \ref{Table:foslp_remedy} presents the results obtained for the FOSLP using different strategies. 
\begin{table*}[htb]
  \centering
  \caption{Quantitative comparison of the results by solving different shooting functions}
  \setlength\tabcolsep{2pt} 
  \begin{tabular}{cccc}
  \hline
  Item      & $\boldsymbol{\Phi}^f_F(\boldsymbol{z}^{f}_{\text{F,ICVN}}(\delta))$ in (\ref{EQ:TPBVP_law_ICVN_fuel})    
  & $\boldsymbol{\Phi}^b_h(\boldsymbol z^b_{\text{h,SICVN}}(\kappa,\delta))$ in (\ref{EQ:TPBVP_new_homo})   
  & $\boldsymbol{\Phi}^f_h(\boldsymbol{z}^{f}_{\text{h,SICVN}}(\kappa,\delta))$ in (\ref{EQ:TPBVP_new_homo_forwaed}) \\ 
  \hline
  No. of successful landings & $8,892$         &$9,952$  &\textcolor{blue}{$9,930$}\\ 
  No. of convergent solutions with $r(t) < R_0$ & $57$   &$0$    &\textcolor{blue}{$0$}\\ 
  Average computational time (s) & $0.7125$   &$0.6218$    &\textcolor{blue}{$0.8348$}\\ 
  Average No. of iterations & $36.63$   &$30.84$   &\textcolor{blue}{$32.32$}\\ 
  Average No. of function evaluations & $262.97$  &$238.72$ &\textcolor{blue}{$241.59$}\\ 
  \textcolor{black}{Success rate ($\%$)}  & $89.35$    & $100$ &\textcolor{blue}{$99.78$}\\
  \hline
  \label{Table:foslp_remedy}
  \end{tabular}
  \end{table*}
A total of $9,952$ feasible fuel-optimal solutions are found using the homotopy approach in (\ref{EQ:TPBVP_new_homo}), indicating a success rate of $100\%$. In contrast, the conventional way only finds $8,892$ feasible solutions, and $57$ convergent solutions are infeasible, indicated by $r(t) < R_0$. This discrepancy is likely due to poor initial guess of the \textcolor{black}{shooting vector} in (\ref{EQ:solutionspace_ICVN_fuel}). Consequently, the conventional indirect method achieves a success rate of $89.35\%$.
Keep in mind that implementing the homotopy approach in (\ref{EQ:TPBVP_new_homo}) requires finding the solution to the TOSLP first, which has an average computational time of $0.0557$ seconds. Therefore, for a feasible initial condition, the homotopy approach in (\ref{EQ:TPBVP_new_homo}) can find the fuel-optimal solution within a total computational time of $0.6775$ seconds, which is even faster than the conventional indirect method ($0.7125$ seconds) \textcolor{black}{despite} an extra homotopy parameter $\kappa$. Remarkably, this computational efficiency is comparable to the learning-based method in \cite{you2021learning}.
\textcolor{blue}{Moreover, we can observe that $9,930$ feasible fuel-optimal solutions are obtained by solving (\ref{EQ:TPBVP_new_homo_forwaed}), which is expected because the $9,930$ feasible time-optimal solutions found by solving (\ref{EQ:TPBVP_law_ICVN_pm_p0}) can provide accurate initialization of the shooting vector for (\ref{EQ:TPBVP_new_homo_forwaed}). 
In summary, solving (\ref{EQ:TPBVP_new_homo_forwaed}), in which the dynamics are propagated forward, requires a total computational time of $0.9303$ seconds to find the solution to the TOSLP. In contrast, our proposed method decreases the average total computational time by $27.17\%$.
}
\section{CONCLUSIONS} \label{Colus}
In this work, we introduced a novel approach called Physics-Informed Indirect Method (PIIM) to solve trajectory optimization problems quickly and robustly. \textcolor{blue}{Its key feature is the incorporation of physics-informed information into all shooting variables. 
An analytical method was provided to estimate the minimum flight time. By eliminating the mass co-state and the numerical factor using a physical fact at the final time, the co-state vector at the final time was constrained on a unit 3-D hypersphere.
By further propagating the dynamics backward, the physical significance of the optimal control was exploited to reduce the solution space. All the physics-informed information was embedded into a shooting function, enabling the Time-Optimal Soft Landing Problem (TOSLP) to be solved quickly and robustly.} Furthermore, the PIIM was applied to a Fuel-Optimal Soft Landing Problem (FOSLP) via a homotopy approach. 

\textcolor{blue}{Numerical simulations demonstrated that
the proposed estimation of the minimum flight time improved the convergence rate, and it reduced the computational time by $56.88\%$ to solve the TOSLP. To ensure that the final time of the convergent solution is positive, we proposed a simple remedy strategy. In this way, the computational time by using the PIIM to solve the TOSLP was reduced to $0.0557$ seconds with a success rate of $100\%$, compared to  $0.2463$ seconds with a success rate of $55.28\%$ by using the conventional indirect method. The results obtained by the homotopy approach showed that the FOSLP could be solved more quickly and robustly, which was identified by a computational time of $0.6775$ seconds with a success rate of $100\%$, compared to $0.7125$ seconds with a success rate of $89.35\%$ by using the conventional indirect method.}

Future research directions include the generalization of the PIIM to time- and fuel-optimal problems, as well as its applicability to facilitating the dataset generation for training neural networks in order to derive the real-time optimal solution.
\section{ACKNOWLEDGMENT}
The authors are grateful to Roberto Armellin, Adam Evans, and the reviewers for their valuable  suggestions on improving the paper.

\bibliography{main_arxiv}

\end{document}